\documentclass[10pt,letterpaper]{article}
\usepackage{arxiv}
\usepackage[T1]{fontenc}
\usepackage{amsmath,amsthm,mathtools}
\usepackage{newtxmath}
\usepackage{booktabs,tabularx,array,multirow}
\usepackage{graphicx,tikz,xcolor}
\usetikzlibrary{arrows.meta,positioning,fit,calc,backgrounds}
\usepackage{microtype,enumitem,siunitx,xspace,url}
\usepackage[numbers,sort&compress]{natbib}
\usepackage[hidelinks]{hyperref}
\usepackage{placeins,flafter,float,needspace,listings}
\usepackage[font=small,labelfont=bf]{caption}
\usepackage{textcomp}
\newcolumntype{Y}{>{\raggedright\arraybackslash}X}
\newcolumntype{P}[1]{>{\raggedright\arraybackslash}p{#1}}
\newcolumntype{C}[1]{>{\centering\arraybackslash}p{#1}}
\newcommand{\R}{\mathbb{R}}

\newcommand{\Z}{\mathbb{Z}}

\newcommand{\optrel}{\simeq_{\mathrm{opt}}}
\newcommand{\dirmap}{\rightsquigarrow}
\newtheorem{theorem}{Theorem}

\newtheorem{proposition}{Proposition}
\newtheorem{reconstructionlemma}{Lemma}[section]

\theoremstyle{definition}
\newtheorem{definition}{Definition}
\definecolor{rpxblue}{HTML}{245B85}
\definecolor{rpxteal}{HTML}{1A786D}
\definecolor{rpxorange}{HTML}{B36B2C}

\newcommand{\RPXKernelTotalSeconds}{472}
\newcommand{\RPXKernelMedianSeconds}{0.288}
\newcommand{\RPXKernelPninetyfiveSeconds}{9.99}
\newcommand{\RPXKernelMaxSeconds}{55.6}
\newcommand{\RPXCheckingTotalSeconds}{510}
\newcommand{\RPXCheckingMedianSeconds}{0.312}
\newcommand{\RPXCheckingPninetyfiveSeconds}{10.1}
\newcommand{\RPXCheckingMaxSeconds}{57.8}
\newcommand{\RPXWallTotalSeconds}{548}
\newcommand{\RPXWallMedianSeconds}{0.382}
\newcommand{\RPXWallPninetyfiveSeconds}{11.6}
\newcommand{\RPXWallMaxSeconds}{61}
\newcommand{\RPXCLITotalTotalSeconds}{544}

\newcommand{\RPXRSSGiBMedian}{0.0673}
\newcommand{\RPXRSSGiBPninetyfive}{1.11}
\newcommand{\RPXRSSGiBMax}{4.16}
\newcommand{\RPXProofThroughputMB}{22.3}

\newcommand{\RPXRSSMiBMedian}{69}
\newcommand{\RPXSquareCheckingSeconds}{10.1}
\newcommand{\RPXNFiveCheckingMilliseconds}{7.95}
\newcommand{\RPXMenuNativeCheckingMilliseconds}{0.074}
\newcommand{\RPXMenuNativeMinMilliseconds}{0.073}
\newcommand{\RPXMenuNativeMaxMilliseconds}{0.076}
\newcommand{\RPXMenuNativeProcessMilliseconds}{3.61}
\newcommand{\RPXMenuNativeCheckingMicroseconds}{74}
\newcommand{\RPXMenuCoreCheckingMilliseconds}{126.4}
\newcommand{\RPXMenuCoreMinMilliseconds}{126.2}
\newcommand{\RPXMenuCoreMaxMilliseconds}{127.1}
\newcommand{\RPXMenuCoreProcessMilliseconds}{130.6}

\newcommand{\RPXTitle}{RPX: Certifying Agent-Generated MILP Reformulations at Scale}
\newcommand{\RPXAuthors}{Haozhen Bai, Wenjie Zhu, Qian Hu, and Caihua Chen}
\hypersetup{pdftitle={\RPXTitle},pdfauthor={\RPXAuthors}}
\renewcommand{\shorttitle}{RPX: Certifying MILP Reformulations}
\title{\RPXTitle}
\author{%
  \RPXAuthors\\[0.35em]
  \normalfont School of Management and Engineering, Nanjing University\\
  \normalfont Nanjing, China\\
  \small
  \href{mailto:baihaozhen@smail.nju.edu.cn}{\texttt{baihaozhen@smail.nju.edu.cn}}\quad
  \href{mailto:ZhuWenjie@smail.nju.edu.cn}{\texttt{ZhuWenjie@smail.nju.edu.cn}}\quad
  \href{mailto:huqian@nju.edu.cn}{\texttt{huqian@nju.edu.cn}}\quad
  \href{mailto:chchen@nju.edu.cn}{\texttt{chchen@nju.edu.cn}}%
}
\date{}
\begin{document}
\maketitle
\begin{abstract}
We present RPX, an independent certificate checker for AI-generated reformulations of mixed-integer
linear programs (MILPs). RPX verifies that every feasible solution of either the original or the
reformulated model can be converted into a feasible solution of the other with no worse objective value.
This guarantees equal optimal values and enables solution recovery in the original variables.
The checker derives the required checks directly from both models and verifies the certificate using
exact arithmetic and reusable reconstruction rules. RPX certifies an agent-generated reformulation
for each of the 240 instances in MIPLIB 2017 Benchmark v2. The source models contain up to 27.3 million
matrix nonzeros. The median time to decode and check a certificate is \RPXCheckingMedianSeconds{}~s,
with a 95th percentile of \RPXCheckingPninetyfiveSeconds{}~s.
\end{abstract}
\section{Introduction}
\label{sec:introduction}

AI agents can now reformulate optimization models directly. Such changes must preserve
the intended optimization problem. This issue is becoming practical as large language models (LLMs) move
beyond describing mathematical programs to producing and modifying solver-ready models.
Systems such as OptiMUS and autoformulation generate optimization models and solver code
from natural language~\cite{ahmaditeshnizi2024optimus,astorga2025autoformulation}, while agents
can also inspect numerical MILP structure, eliminate variables, aggregate constraints, and produce
new formulations. More broadly, LLM4MIP reports resolving 52 of 217 open MIPLIB instances using
LLM-assisted solver workflows~\cite{huang2026llm4mip}. As agents gain greater freedom to modify executable
optimization artifacts, independently establishing the correctness of their reformulations becomes increasingly important.

The difficulty is that an incorrect reformulation can look successful. A missing constraint or an invalid integrality
relaxation may make the target easier to solve by changing its feasible set; a lost objective constant can instead
distort comparisons of objective values. An agent's explanation does not provide a machine-checkable guarantee,
and agreement between numerical solver runs tests only particular outcomes.
At the same time, requiring the two models to have identical feasible sets is too restrictive: useful reformulations
may intentionally remove dominated solutions or introduce a different set of variables. What is needed is therefore
a precise correctness criterion that is strong enough to preserve the optimization result, flexible enough to permit
such reformulations, and independently checkable from the resulting models.

\paragraph{The missing verification interface.} Existing work studies reformulation verification through mapped-solution
tests, formal theorem proving, and
SMT~\cite{zhai2025equivamap,robbins2026flare,bhattacharyya2026sover};
proof logging also certifies solver reasoning and presolve
reductions~\cite{cheung2017vipr,hoen2024presolve}.
For open-ended agent reformulation, we seek an interface that takes the
final source and target MILPs together with a certificate, without
requiring the agent's search history or solving either model.
The certificate must establish correctness for every feasible solution,
including how solutions can be reconstructed when variables are
eliminated or replaced.

This interface must also scale: our evaluation includes source models
with up to 1.43 million variables and 27.3 million matrix nonzeros.
Certification must therefore support changes in variable representation
while checking the complete models exactly and efficiently.

\paragraph{Our approach.}

We introduce RPX, an independent checker for this interface.
Its fixed criterion, \texttt{MI\_OPT}, requires that every feasible
solution of either model can be converted into a feasible solution
of the other with no worse objective value.
This preserves feasibility and the optimal objective value (Section~2).
RPX generates checks for every variable bound, integrality condition,
constraint, and objective comparison from the models themselves.

\begin{figure}[htbp]
\centering
\begin{tikzpicture}[font=\small,arr/.style={-{Latex[length=2mm]},thick},
 box/.style={draw,rounded corners=2pt,align=center,minimum height=10mm,text width=26mm,inner sep=3pt},
 model/.style={box,fill=rpxblue!7},
 untrusted/.style={box,draw=rpxorange,dashed,fill=rpxorange!7}]
\node[model] (source) at (0,0) {Source MILP};
\node[untrusted,text width=29mm] (agent) at (6.75,0) {AI agent};
\node[model,text width=29mm] (target) at (13.5,0) {Reformulated\\MILP};
\node[font=\scriptsize,text=rpxorange,align=center] at (6.75,0.95)
  {\textbf{Untrusted:} reformulation agent + certificate generation};
\draw[arr] (source)--(agent);
\draw[arr] (agent)--(target);

\node[model,text width=22mm] (source-input) at (0,-2.0) {Source MILP};
\node[font=\large] at (1.5,-2.0) {$+$};
\node[model,text width=29mm] (target-input) at (3.2,-2.0) {Reformulated\\MILP};
\node[font=\large] at (5.0,-2.0) {$+$};
\node[untrusted,text width=23mm] (certificate) at (6.5,-2.0) {Certificate};
\node[font=\scriptsize,text=rpxorange,align=center,text width=45mm] at (6.5,-3.1)
  {Certificate may be produced\\by any method};
\node[box,fill=rpxteal!12,draw=rpxteal,very thick,text width=25mm] (check) at (10.2,-2.0) {RPX checker};
\node[font=\scriptsize,text=rpxteal,align=center,text width=34mm,inner sep=0pt] (trust-label) at (10.2,-1.2)
  {\textbf{Trusted:} exact model\\interpretation + RPX checker};
\node[box,text width=22mm] (out) at (13.5,-2.0) {Accepted\\Invalid\\ResourceLimit};
\draw[arr] (certificate)--(check);
\draw[arr] (check)--(out);
\begin{scope}[on background layer]
\node[draw=rpxteal,dashed,rounded corners,fill=rpxteal!3,
  fit=(trust-label)(check),inner sep=2mm] {};
\end{scope}
\end{tikzpicture}
\caption{\textbf{Flexible generation, independent acceptance.} Given the complete source MILP, reformulated MILP, and certificate, RPX returns Accepted, Invalid, or ResourceLimit. Orange dashed boxes denote untrusted generation; the teal boundary denotes trusted exact model interpretation and checking.}
\label{fig:architecture}
\end{figure}
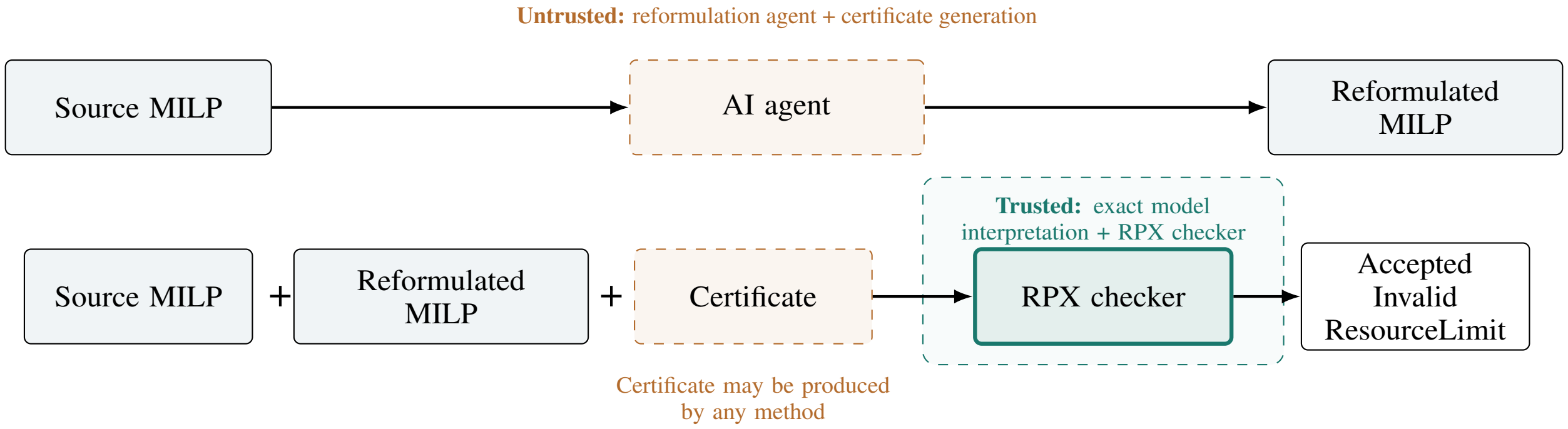

A small exact Core checks linear arithmetic, integer reasoning,
and case analysis.
Fixed reconstruction rules avoid enumerating combinations of discrete
choices and support joint recovery of variables linked by shared
constraints.
This design separates discovery from acceptance: agents may search
freely, but only the independent checker determines whether a
reformulation is certified (Figure~\ref{fig:architecture}).

Our contributions are:
\begin{enumerate}[leftmargin=*,itemsep=2pt,topsep=3pt]
\item \textbf{A fixed certification interface for open-ended agent reformulation.} We formulate correctness as a
property of the final source–target model pair. The \texttt{MI\_OPT} contract
permits useful changes such as removing dominated solutions and changing the variable space, while preserving
feasibility and the optimal objective value.
\item \textbf{An exact and scalable checker for this interface.} RPX derives the required checks from the complete
models and combines a small exact reasoning core with fixed reconstruction rules for reformulations that eliminate
or reorganize variables.
\item \textbf{Evidence that final-pair certification is practical at MIPLIB scale.} We evaluate RPX on 240 source–target
pairs originating from agent reformulation runs over all MIPLIB 2017 Benchmark v2 instances~\cite{gleixner2021miplib}. The same checker and fixed
rule set certify all 240 pairs.
Median decoding-and-checking time is \RPXCheckingMedianSeconds{}~s and the 95th percentile is \RPXCheckingPninetyfiveSeconds{}~s. A separate robustness campaign rejects
all 12,544 confirmed-invalid certificates.
\end{enumerate}

\section{The \texttt{MI\_OPT} Correctness Criterion}
\label{sec:semantics}

RPX needs a fixed criterion for deciding whether a reformulated model is correct. Requiring the source and target
models to have identical feasible sets, or even a bijection between their feasible solutions, would be too restrictive:
a useful reformulation may deliberately remove solutions that can never improve the objective. We first illustrate this
point with a small example, and then formalize the criterion used by RPX.

\subsection{Why identical feasible sets are too restrictive}
\label{sec:example}
Consider a simple network in which an integer quantity $q\in\{0,1,2\}$ must pass through two consecutive links.
Opening the first link costs 1 and opening the second costs 2; each has capacity 2, and each unit delivered earns 2.
With binary opening decisions $a,b$, the source is
\begin{equation}
\min\ a+2b-2q\quad\text{s.t. }q\le2a,\ q\le2b.
\label{eq:simple-source}
\end{equation}
An agent replaces the two decisions by one binary $z$:
\begin{equation}
\min\ 3z-2q\quad\text{s.t. }q\le2z.
\label{eq:simple-target}
\end{equation}
Both models retain $q\in\{0,1,2\}$.

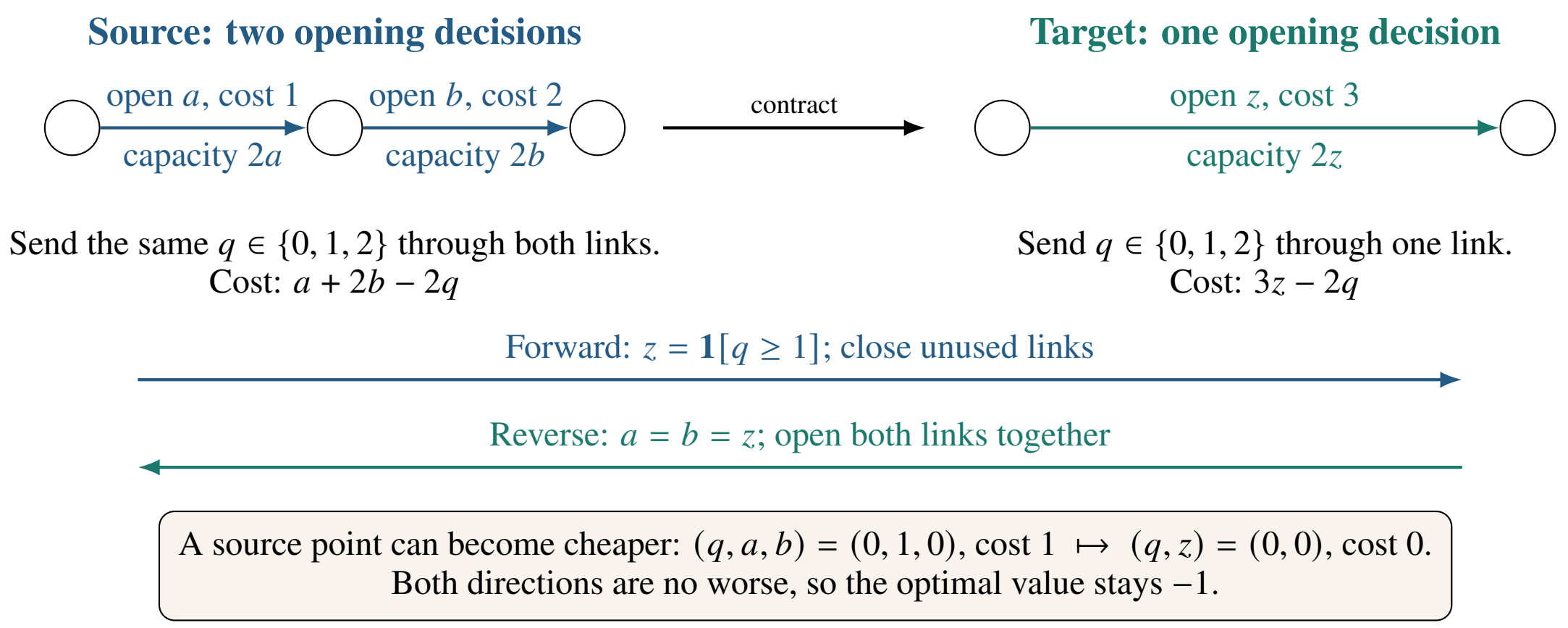
\begin{figure}[htbp]
\centering
\begin{tikzpicture}[font=\small,arr/.style={-{Latex[length=2mm]},thick},
 vertex/.style={circle,draw,fill=white,minimum size=5mm,inner sep=0pt}]
\node[font=\bfseries,text=rpxblue] at (2.8,1.4) {Source: two opening decisions};
\node[font=\bfseries,text=rpxteal] at (11.3,1.4) {Target: one opening decision};
\node[vertex] (s) at (0.4,0.55) {};
\node[vertex] (m) at (2.8,0.55) {};
\node[vertex] (t) at (5.2,0.55) {};
\draw[arr,rpxblue] (s)--node[above]{open $a$, cost $1$}node[below]{capacity $2a$}(m);
\draw[arr,rpxblue] (m)--node[above]{open $b$, cost $2$}node[below]{capacity $2b$}(t);
\node[vertex] (ss) at (8.9,0.55) {};
\node[vertex] (tt) at (13.7,0.55) {};
\draw[arr,rpxteal] (ss)--node[above]{open $z$, cost $3$}node[below]{capacity $2z$}(tt);
\draw[arr] (5.8,0.55)--node[above,font=\scriptsize]{contract}(8.2,0.55);
\node[align=center] at (2.8,-0.7) {Send the same $q\in\{0,1,2\}$ through both links.\\Cost: $a+2b-2q$};
\node[align=center] at (11.3,-0.7) {Send $q\in\{0,1,2\}$ through one link.\\Cost: $3z-2q$};
\draw[arr,rpxblue] (1.0,-1.75)--node[above]{Forward: $z=\mathbf{1}[q\ge1]$; close unused links}(13.1,-1.75);
\draw[arr,rpxteal] (13.1,-2.55)--node[above]{Reverse: $a=b=z$; open both links together}(1.0,-2.55);
\node[draw,rounded corners,fill=rpxorange!8,align=center,inner sep=5pt] at (7.1,-3.45)
{A source point can become cheaper: $(q,a,b)=(0,1,0)$, cost $1$ $\;\mapsto\;$ $(q,z)=(0,0)$, cost $0$.\\Both directions are no worse, so the optimal value stays $-1$.};
\end{tikzpicture}
\caption{\textbf{Preserve the optimization task while discarding unnecessary decisions.} The forward map removes wasteful activations; the reverse map reconstructs a valid original solution. No bijection is needed. This tiny example captures the same recovery obligation that large reformulations must satisfy.}
\label{fig:example}
\end{figure}

\textbf{Source to target.} Keep $q$ and set $z=0$ when $q=0$, and $z=1$ otherwise.
If $q=0$, closing both links cannot increase cost.
If $q\ge1$, feasibility and binary domains force $a=b=1$, so the target has the same cost.
\textbf{Target to source.} Set $a=b=z$ and keep $q$; all constraints hold and the cost is unchanged.
Figure~\ref{fig:example} visualizes this proof.

The forward construction discards unnecessary openings: $(q,a,b)=(0,1,0)$ maps to $(q,z)=(0,0)$ and becomes cheaper.
The two models therefore do not have identical feasible sets under a one-to-one correspondence. Nevertheless,
every feasible source solution has a target solution with no greater cost, and every feasible target solution has
a source solution with no greater cost. These two properties are enough to force the same optimal objective value.

\subsection{The \texttt{MI\_OPT} guarantee}
We now formalize the property illustrated above. A rational mixed-integer linear program (MILP) $M$ has feasible set
\[
X_M=\{x:\ell\le x\le u,\ b^-\le Ax\le b^+,\ x_I\in\Z^I\}
\]
and minimization objective $f_M(x)=c^\top x+c_0$.
Finite coefficients are interpreted as exact rationals; continuous coordinates range over $\R$.
Absent bounds denote infinity, and maximization is normalized by negating the complete objective, including its constant.
Thus a decimal such as $0.333$ denotes $333/1000$.

\begin{definition}[No-worse construction]
\label{def:domination}
We write $M\dirmap N$ when there exists a total map $F:X_M\to X_N$ with
\[
f_N(F(x))\le f_M(x)\qquad\text{for all }x\in X_M.
\]
The fixed goal $M\optrel N$ requires $M\dirmap N$ and $N\dirmap M$.
\end{definition}

\begin{theorem}[Optimization preservation]
\label{thm:equivalence}
If $M\optrel N$, the models are simultaneously feasible or infeasible, and their extended optimal values are equal:
\[
\inf_{x\in X_M}f_M(x)=\inf_{y\in X_N}f_N(y),
\]
where $\inf\varnothing=+\infty$.
Objective unboundedness is preserved, and either construction sends an attained optimum to an optimum of the other model.
\end{theorem}
\begin{proof}
A total map sends any feasible input to a feasible output, giving feasibility in both directions.
The forward inequality gives $\inf f_N\le\inf f_M$; the reverse gives the opposite inequality, also when the value is $-\infty$.
For an attained optimum, a feasible no-worse image must attain this common value.
\end{proof}

\texttt{MI\_OPT} therefore preserves the optimization outcome without requiring identical variable representations, feasible sets, or LP relaxations.

\section{How RPX Certifies \texttt{MI\_OPT}}
\label{sec:method}

RPX checks a construction in each direction required by \texttt{MI\_OPT}. We describe the source-to-target direction;
the reverse direction exchanges the two models.

\subsection{Checking a directional construction}
\label{sec:check-construction}
A directional certificate describes how to construct the target variables from an arbitrary feasible source solution.
The source variables remain symbolic, so the construction must work for every feasible source solution.

\begin{figure}[htbp]
\centering
\begin{tikzpicture}[font=\small,
 arr/.style={-{Latex[length=1.7mm]},thick,shorten <=1.5pt,shorten >=1.5pt},
 box/.style={draw,rounded corners=2pt,align=center,inner xsep=5pt,inner ysep=5pt},
 branch/.style={box,text width=60mm,minimum height=17mm,fill=rpxteal!6}]
\node[box,text width=31mm,fill=rpxblue!8] (source) at (0,0)
  {\textbf{Source facts}\\[2pt]
   $q\in\{0,1,2\}$\\
   $a,b\in\{0,1\}$\\
   $q\le 2a,\ q\le 2b$};
\node[branch] (zero) at (6.0,1.04)
  {\textbf{Case $q=0$}\\[2pt]
   Set $z=0$. Check $q\le 2z$, $z$ binary,\\
   and $3z-2q\le a+2b-2q$.};
\node[branch] (positive) at (6.0,-1.04)
  {\textbf{Case $q\ge1$}\\[2pt]
   From $2a,2b\ge1$, round to $a=b=1$.\\
   Set $z=1$; check capacity, domain, cost.};
\node[box,text width=34mm,fill=rpxteal!12] (result) at (12.2,0)
  {\textbf{Forward guarantee}\\[2pt]
   Both cases satisfy the output domains, rows, and objective inequality.};
\draw[arr] (source.east)--(zero.west);
\draw[arr] (source.east)--(positive.west);
\draw[arr] (zero.east)--(result.west);
\draw[arr] (positive.east)--(result.west);
\end{tikzpicture}
\caption{\textbf{Small checks compose into a full guarantee.} Forward certificate for the running example.
Both branches are checked, with no sharing of branch-local assumptions. A separate reverse proof is also required.}
\label{fig:proof-dag}
\end{figure}
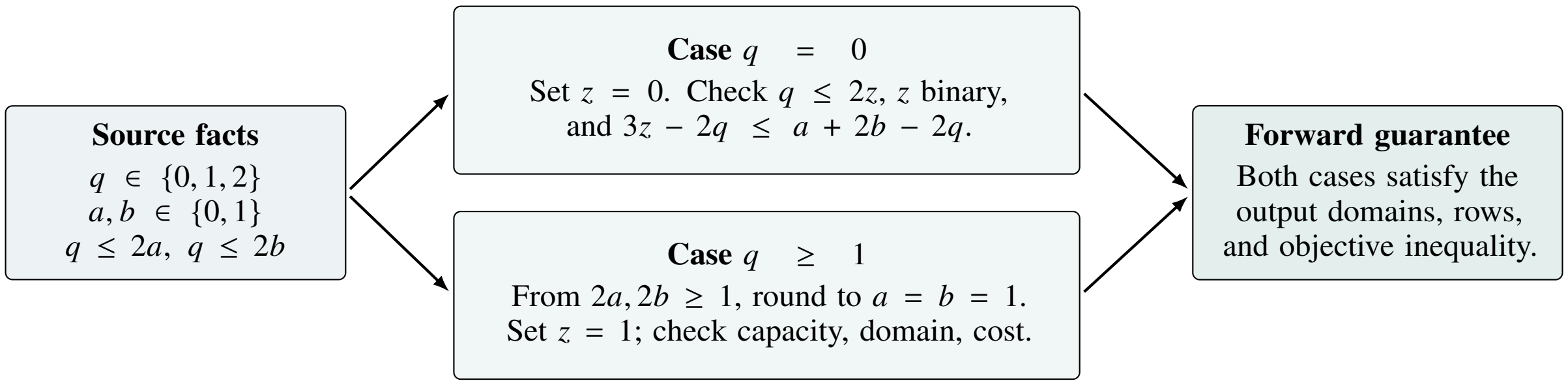

Figure~\ref{fig:proof-dag} illustrates the checks for the running example. The forward construction sets \(z=0\) when \(q=0\) and \(z=1\)
when \(q\ge1\). In each case, RPX establishes that \(z\) is binary, the target capacity constraint holds, and the objective
does not increase. The second case uses integer reasoning to derive \(a=b=1\) from the source constraints.

More generally, RPX generates every required variable bound, integrality condition, constraint, and objective comparison
from the target model and checks them for the constructed variables. The certificate must justify this complete set of checks.

\subsection{Core reasoning}
\label{sec:core}

The Core provides the elementary exact reasoning used to establish the requirements in Section~\ref{sec:check-construction}.
It works with affine expressions, meaning linear combinations of variables plus a constant. Specifically, it records three kinds of statements:
an inequality \(e\le 0\), an equality \(e=0\), or the fact that \(e\) is integer-valued. A certificate may combine valid inequalities and
equalities using exact rational arithmetic, use integrality to round a bound, or justify a weaker statement from known variable bounds.
For example, \(2a\ge 1\) together with \(a\in\mathbb{Z}\) implies \(a\ge 1\).
The Weaken instruction packages a variable-bound implication into one check, avoiding explicit multipliers for every participating bound.

The Core also supports exhaustive case analysis. Different cases may use different constructions, but assumptions
introduced within a case remain local to that case. Appendix~\ref{app:core-details} specifies the supported splits and their acceptance conditions.

These rules support common reasoning such as substitution, aggregation, integer rounding, and finite case analysis without introducing a
separate rule for every named reformulation technique. Appendix~\ref{app:core} gives their exact acceptance conditions and a formal characterization
of the finite case-based constructions they can certify.

\subsection{Reusable reconstruction rules}
\label{sec:lifts}

Elementary case analysis becomes cumbersome when repeated discrete selections duplicate subsequent checks,
or when several variables must be recovered jointly to satisfy shared constraints. RPX provides six fixed
reconstruction rules for these tasks (Table~\ref{tab:lifts-main}).

The certificate establishes a rule's preconditions using the Core, then uses the reconstructed variables
and the relations guaranteed by that rule to establish the complete model
checks in Section~\ref{sec:check-construction}. Appendix~\ref{app:lifts} gives the exact rule specifications.

\begin{table}[htbp]
\caption{Fixed reconstruction rules. The first three avoid repeated case splitting for local
choices; the last three reconstruct flow variables that must satisfy shared network constraints.}
\label{tab:lifts-main}
\centering\small
\begin{tabularx}{\linewidth}{@{}lY@{}}
\toprule
Rule & What it reconstructs\\
\midrule
Interval & One continuous variable satisfying several established lower and upper bounds.\\
Menu & Binary option variables from an on/off indicator and the index of the selected option.\\
Signed & A binary direction indicator consistent with the sign of a bounded quantity.\\
\midrule
Flow & Destination-specific flows whose sum equals the aggregate flow on every edge.\\
Subflow & A smaller flow that preserves all deliveries and uses no more of any edge.\\
UnitFlowRound & Binary edge decisions that preserve net flow at every node and do not increase cost.\\
\bottomrule
\end{tabularx}
\end{table}

\paragraph{Local example: Menu.}

Suppose the source model uses binary variables \(q_0,\ldots,q_K\) to represent an optional
selection: either exactly one option is selected or no option is selected. The reformulated
target stores only a binary variable \(A\), indicating whether a selection is
active, and an integer \(I\), indicating the selected option.

\begin{figure}[htbp]
\centering
\begin{tikzpicture}[font=\small,
panel/.style={draw=rpxteal!45!black,rounded corners=3pt,minimum width=33mm,
              minimum height=23mm,fill=rpxteal!4},
cell/.style={draw=rpxteal!50!black,fill=white,minimum size=5.3mm,inner sep=0pt},
selected/.style={fill=rpxteal!27,very thick,font=\small\bfseries}]
\foreach \xx/\heading/\aa/\ii/\qa/\qb/\qc in {
  0/{Select 0}/1/0/1/0/0,
  3.85/{Select 1}/1/1/0/1/0,
  7.70/{Select 2}/1/2/0/0/1,
  11.55/{No selection}/0/0/0/0/0}{
  \begin{scope}[shift={(\xx,0)}]
    \ifnum\aa=0
      \node[panel,draw=black!40,fill=black!3] at (0,0) {};
    \else
      \node[panel] at (0,0) {};
    \fi
    \node[font=\small\bfseries] at (0,.77) {\heading};
    \node at (0,.34) {$A=\aa,\ I=\ii$};
    \foreach \dx/\jj/\bit in {-.66/0/\qa,0/1/\qb,.66/2/\qc}{
      \node[font=\scriptsize] at (\dx,-.12) {$q_{\jj}$};
      \ifnum\bit=1
        \node[cell,selected] at (\dx,-.62) {1};
      \else
        \node[cell] at (\dx,-.62) {0};
      \fi
    }
  \end{scope}
}
\end{tikzpicture}
\caption{\textbf{Recover a discrete choice once, then reuse its laws.} Menu reconstructs the optional one-hot vector from activation $A$ and index $I$. For $K=2$, all four feasible input states are shown, including the distinction between selecting index 0 and making no selection.}
\label{fig:menu-example}
\end{figure}
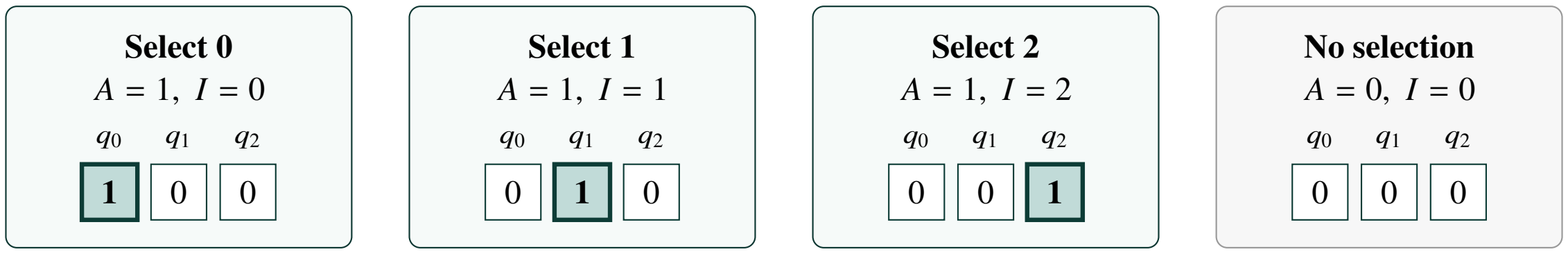

Menu checks that \(A\in\{0,1\}\), that \(I\) is integer-valued, and that \(0\le I\le KA\).
It then reconstructs binary variables \(q_0,\ldots,q_K\) satisfying
\[
\sum_{j=0}^{K}q_j=A,
\qquad
\sum_{j=0}^{K}j q_j=I.
\]
If \(A=1\), exactly option \(I\) is selected. If \(A=0\), all \(q_j\) are zero. The pair \((A,I)\)
therefore distinguishes selecting option \(0\) from selecting no option. Figure~\ref{fig:menu-example} illustrates these possibilities.

Menu records these relations once, allowing subsequent checks to reuse them without explicitly enumerating the possible \((A,I)\) states.

\paragraph{Joint example: UnitFlowRound.}

\begin{figure}[htbp]
\centering
\begin{tikzpicture}[x=1.45cm,font=\small,vertex/.style={circle,draw,fill=white,inner sep=1.5pt,minimum size=5mm},
edge/.style={-{Latex[length=1.5mm]},thick}]
\foreach \shift/\prefix in {0/f,4.0/g}{
  \node[vertex] (\prefix s) at (\shift,0) {$s$};
  \node[vertex] (\prefix u) at (\shift+1,0.8) {$u$};
  \node[vertex] (\prefix v) at (\shift+1,-0.8) {$v$};
  \node[vertex] (\prefix t) at (\shift+2,0) {$t$};
}
\draw[edge,rpxblue] (fs)--node[above left,inner sep=1pt] {$\frac12$}(fu);
\draw[edge,rpxblue] (fu)--node[above right,inner sep=1pt] {$\frac12$}(ft);
\draw[edge,rpxblue] (fs)--node[below left,inner sep=1pt] {$\frac12$}(fv);
\draw[edge,rpxblue] (fv)--node[below right,inner sep=1pt] {$\frac12$}(ft);
\draw[edge,rpxteal,very thick] (gs)--node[above left,inner sep=1pt] {1}(gu);
\draw[edge,rpxteal,very thick] (gu)--node[above right,inner sep=1pt] {1}(gt);
\draw[edge,gray,dashed] (gs)--node[below left,inner sep=1pt] {0}(gv);
\draw[edge,gray,dashed] (gv)--node[below right,inner sep=1pt] {0}(gt);
\draw[edge] (2.5,0)--(3.5,0);
\node[font=\small] at (3,0.4) {Joint choice};
\node[font=\scriptsize] at (1,1.25) {upper path: cost 2};
\node[font=\scriptsize] at (1,-1.2) {lower path: cost 4};
\node[align=center] at (1,-1.9) {Fractional $f$\\cost $3$};
\node[align=center] at (5,-1.9) {Binary $g$\\cost $2$};
\end{tikzpicture}
\caption{\textbf{Why recovery can require a joint construction.} In this schematic unit-capacity network, the upper path costs
2 and the lower path costs 4. Half a unit on each path costs 3; choosing the whole upper path costs 2. Both send one
unit from $s$ to $t$ and preserve the zero balance at $u,v$. Independent arc rounding can break these balances.
UnitFlowRound certifies existence of a binary, no-worse flow after checking its premises.}
\label{fig:joint-recovery}
\end{figure}
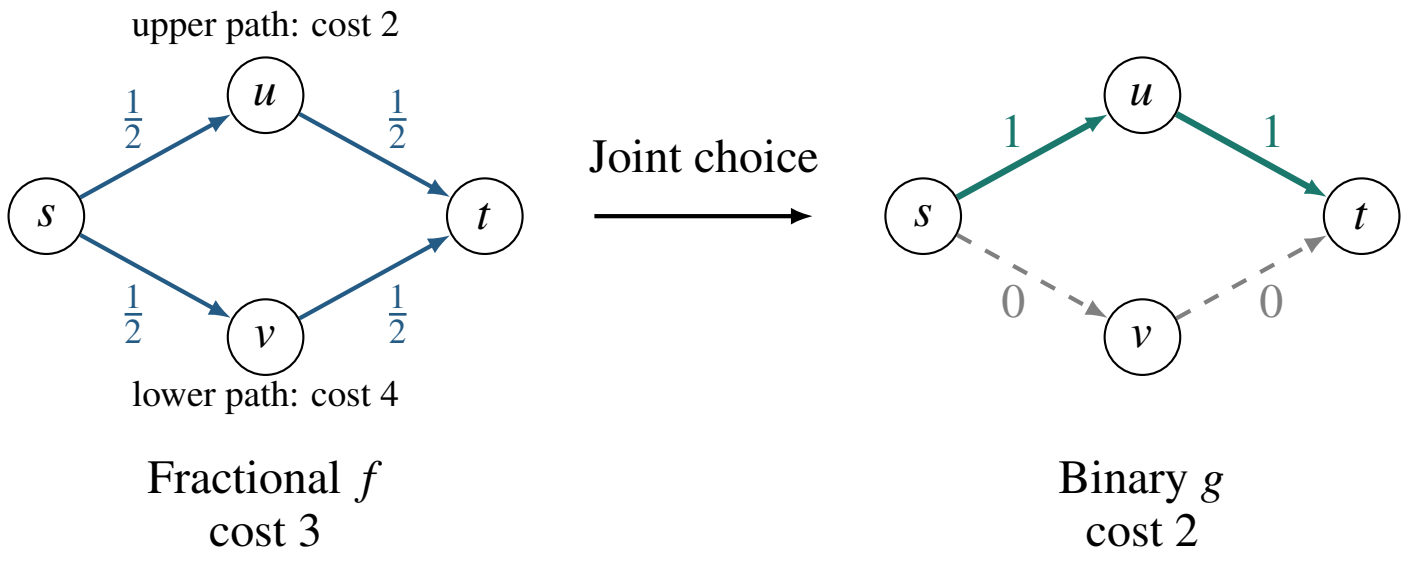

Suppose a reformulated model contains edge variables \(f\in[0,1]^E\), while
the source model requires binary edge decisions. Choosing each edge independently may violate flow conservation.
UnitFlowRound instead reconstructs all edge decisions together. When the net flow at every node is integer-valued,
it supplies binary variables \(g\in\{0,1\}^E\) with the same net flow at every node and no greater total cost.
Figure~\ref{fig:joint-recovery} illustrates why this joint choice is necessary.

\FloatBarrier

\section{Soundness and the Acceptance Boundary}
\label{sec:implementation}

We now state what acceptance guarantees and identify the components on which that guarantee depends.

\paragraph{Acceptance boundary.} RPX checks two models represented with exact rational coefficients and one binary
certificate. It performs no optimization and no proof search. The checker returns one of three outcomes:
Accepted certifies \texttt{MI\_OPT}, Invalid reports a failed certificate check, and ResourceLimit reports an exhausted checking budget.

Table~\ref{tab:acceptance-boundary} identifies the trusted components. When checking starts from MPS or LP text, exact import is also trusted.
Appendix~\ref{app:robustness} reports implementation validation.

\begin{table}[htbp]
\caption{The RPX acceptance boundary. Components inside the boundary must correctly implement the specified checks; components outside it may use arbitrary methods.}
\label{tab:acceptance-boundary}
\centering\small
\begin{tabularx}{\textwidth}{@{}YY@{}}
\toprule
\textbf{Inside the acceptance boundary} & \textbf{Outside the acceptance boundary}\\
\midrule
Exact model interpretation & Reformulation agent\\
RPX decoding & Certificate author\\
Exact rational arithmetic & Python authoring API\\
Core and reconstruction-rule implementations & RPA compiler\\
Case scoping & Solvers or scripts used during discovery\\
Automatic generation of all target-model checks & Natural-language explanations\\
\bottomrule
\end{tabularx}
\end{table}

\paragraph{Soundness.}
\label{sec:soundness-main}
\begin{theorem}[RPX soundness]
\label{thm:certificate-soundness}
Assume that exact model interpretation and the RPX acceptance boundary correctly implement the formal rules
specified in Appendices~\ref{app:core-details} and \ref{app:lift-details}. If RPX accepts a certificate for models \(M\) and \(N\), then
\[
M \simeq_{\mathrm{opt}} N.
\]
\end{theorem}
\begin{proof}
Consider one direction, from \(M\) to \(N\). Each Core step preserves a statement that holds for every feasible input solution.
Case splits cover all feasible inputs, and assumptions introduced in one case cannot be used in another. Each reconstruction
rule supplies variables satisfying its fixed guarantee once its checked preconditions hold. At the end of every feasible case,
the checks generated from \(N\) establish that the constructed solution is feasible and has no larger objective value.
The accepted directional certificate therefore establishes \(M\rightsquigarrow N\). The reverse direction similarly establishes
\(N\rightsquigarrow M\), and Theorem 1 gives \(M\simeq_{\mathrm{opt}}N\). A complete directional proof appears in Appendix~\ref{app:directional-soundness}.
\end{proof}

\paragraph{Certificate production.} Authors may use a Python interface and the readable RPA representation to produce binary RPX certificates.
RPX independently checks the compiler's output against both models, placing the compiler outside the acceptance boundary.
Appendix~\ref{app:worked-certificates} gives worked certificates; Appendix~\ref{app:robustness} describes
the reproduction workflow.

\section{Evaluation}
\label{sec:evaluation}

We evaluate RPX on 240 model pairs and certificates, measuring coverage, checking cost,
the benefit of reconstruction rules, and enforcement of the acceptance boundary.
Appendix~\ref{app:evaluation-details} gives the experimental setup and supplementary results.

\subsection{Benchmark set and measurement protocol}
The evaluation set contains \(240\) source-target pairs, one for each instance in the MIPLIB 2017 Benchmark v2 collection~\cite{gleixner2021miplib}.
The target models were produced by an AI reformulation agent, and certificates were prepared after target selection.
Source models range from \(30\) to \(1,429,098\) variables, from \(4\) to \(2,897,380\) constraints,
and from \(123\) to \(27,329,856\) matrix nonzeros.

Coverage and timing use one fixed checker configuration for all \(240\) pairs. Timing measurements start from canonical
exact models and prepared binary certificates, and therefore include certificate decoding and checking but exclude
text-model import and certificate production. On an Apple M5 with 16 GiB of memory, each pair is checked in five fresh
processes after one warmup; we report the median for each pair. Appendix~\ref{app:evaluation-details} gives the measurement protocol and complete descriptive results; Appendix~\ref{app:robustness} describes implementation validation and artifact reproduction.

\subsection{Coverage across reformulation patterns}

\begin{figure}[htbp]
\centering
\includegraphics[width=\linewidth]{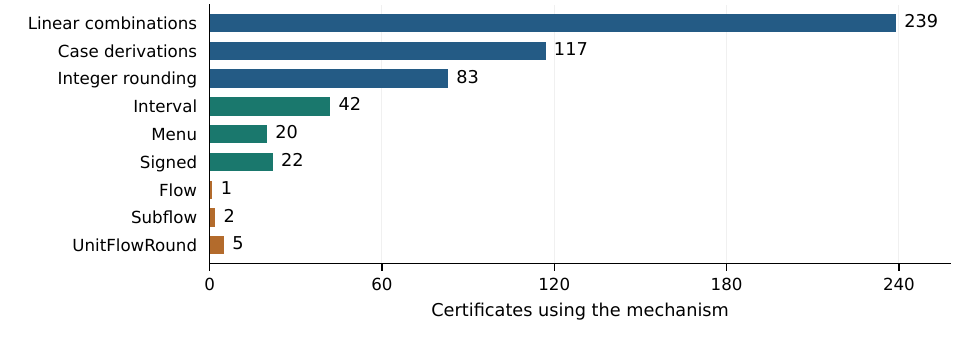}
\caption{\textbf{A shared reasoning vocabulary across reformulations.} Bars count certificates using each mechanism among the 240 certified model pairs. Categories overlap because one certificate can use multiple mechanisms.}
\label{fig:mechanism-usage}
\end{figure}

The same Core and six reconstruction rules certify all 240 supplied pairs: 176 certificates use only the Core and 64 invoke
reconstruction rules. Figure~\ref{fig:mechanism-usage} reports mechanism usage, and Table~\ref{tab:strategies} links
representative reformulations to the mechanisms used to certify them. Median target/source ratios are \(0.886\) for variables, \(0.822\) for constraints, and \(0.748\) for nonzeros,
although some targets grow. Appendix~\ref{appD:statistics} gives the full distributions.

\begin{table}[htbp]
\caption{Representative certified reformulations. The examples illustrate distinct changes in model structure and the RPX mechanisms used to certify them.}
\label{tab:strategies}
\centering\small
\begin{tabularx}{\textwidth}{@{}P{3.1cm}P{3.1cm}YY@{}}
\toprule
Strategy & Example & What the agent changes & What RPX checks\\
\midrule
Algebraic elimination and sparsification & \texttt{square47}, \texttt{square41} & Replaces variables by their known values or formulas, and cancels repeated terms between equations. & Exact linear combinations and complete recovery obligations.\\[3pt]
Boolean structure compression & \texttt{k1mushroom} & Merges variables for identical yes/no conditions and reuses common parts of those conditions. & Integer case reasoning and objective reconstruction.\\[3pt]
Auxiliary-variable relaxation & \texttt{roi5alpha10n8} & Allows fractional auxiliary values and replaces pairwise conflicts with group constraints allowing at most one selection. & Rounding, cases, and Interval witnesses for auxiliary recovery.\\[3pt]
Shipment aggregation & \texttt{n5-3} & Combines shipments from the same source, dropping the separate variables for each destination. & Forward summation and joint Flow decomposition in reverse.\\[3pt]
One-hot block reformulation & \texttt{neos-4532248-waihi} & Allows fractional choices and replaces separate cost constraints with one bound on their weighted sum. & Joint binary recovery with equal balances and no greater cost via UnitFlowRound.\\[3pt]
Global fixing or contradiction & \texttt{enlight\_hard}, \texttt{neos859080} & Fixes all values when only one feasible solution exists, or writes an explicit contradiction when none exists. & Modular linear reasoning, rounding, and exhaustive cases.\\
\bottomrule
\end{tabularx}
\end{table}

Two cases quantify these changes. In \texttt{square47}, algebraic elimination reduces matrix nonzeros from
\(27.3\) million to \(138,331\). The certificate justifies exact cancellations using Core linear reasoning, and RPX
checks the complete pair in \(\RPXSquareCheckingSeconds\)~s. In \texttt{n5-3}, aggregating destination-specific shipments reduces variables
from \(2,550\) to \(850\). The reverse certificate uses the Flow rule to reconstruct shipments preserving every
destination demand and network balance; the complete check takes \(\RPXNFiveCheckingMilliseconds\)~ms.

\FloatBarrier
\subsection{Cost of complete model-pair checking}
\label{sec:checking-cost}

RPX checks every generated variable bound, integrality condition, constraint side, and objective comparison without
solving either MILP. Across the \(240\) pairs, the median certificate decoding-and-checking time is \(\RPXCheckingMedianSeconds\)~s, and \(95\%\) of pairs
complete within \(\RPXCheckingPninetyfiveSeconds\)~s. The maximum is \(\RPXCheckingMaxSeconds\)~s. Median peak memory is \(\RPXRSSMiBMedian\)~MiB, while the largest case reaches \(\RPXRSSGiBMax\)~GiB.

\begin{figure}[htbp]
\centering
\includegraphics[width=\linewidth]{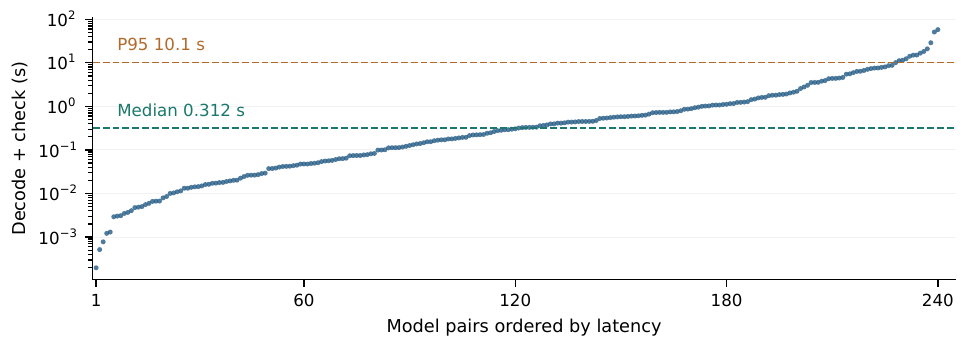}
\caption{\textbf{Distribution of decoding-and-checking time across 240 model pairs.} Points are per-pair medians of five fresh-process
runs after one warmup, ordered by latency. Coverage and timing use the common configuration in Appendix~\ref{appD:protocol}; text-model
import and certificate production are excluded.}
\label{fig:checking-latency}
\end{figure}

\begin{table}[htbp]
\caption{Checking cost over 240 cases. Times are per-case medians of five processes; RSS is each case's largest measured peak. MB is decimal and GiB is binary.}
\label{tab:verification240}
\centering\small
\begin{tabular*}{\linewidth}{@{\extracolsep{\fill}}lrrr@{}}
\toprule
Metric & Median & P95 & Max.\\
\midrule
Certificate (MB) & 8.46 & 233 & 1,059\\
Kernel (s) & \RPXKernelMedianSeconds & \RPXKernelPninetyfiveSeconds & \RPXKernelMaxSeconds\\
Decode + check (s) & \RPXCheckingMedianSeconds & \RPXCheckingPninetyfiveSeconds & \RPXCheckingMaxSeconds\\
Process wall (s) & \RPXWallMedianSeconds & \RPXWallPninetyfiveSeconds & \RPXWallMaxSeconds\\
Peak RSS (GiB) & \RPXRSSGiBMedian & \RPXRSSGiBPninetyfive & \RPXRSSGiBMax\\
\bottomrule
\end{tabular*}
\end{table}

The certificates contain \(57.9\) million checked derivation steps and occupy \(11.397\) GB in total. Summed over
the per-pair medians, decoding and checking takes \(\RPXCheckingTotalSeconds\)~s, corresponding to an observed certificate throughput
of \(\RPXProofThroughputMB\)~MB/s. Figure~\ref{fig:checking-latency} shows the full latency distribution, and Table~\ref{tab:verification240}
reports certificate size, checking time, and memory use. The largest certificate occupies \(1.06\) GB, leaving substantial scope for reducing certificate storage.

\FloatBarrier
\subsection{Why reusable reconstruction rules matter}

To measure the benefit of Menu, we compare it with explicit Core cases encoding the same reconstruction.

The controlled family contains \(k=1,\ldots,6\) independent optional selections. Both certificate variants use the
same model pairs and preserve the same source variables; they differ only in whether each selection is represented
by one Menu invocation or expanded into explicit cases. All certificates are checked by the same production checker
under the same limits.

\begin{figure}[htbp]
\centering
\includegraphics[width=\linewidth]{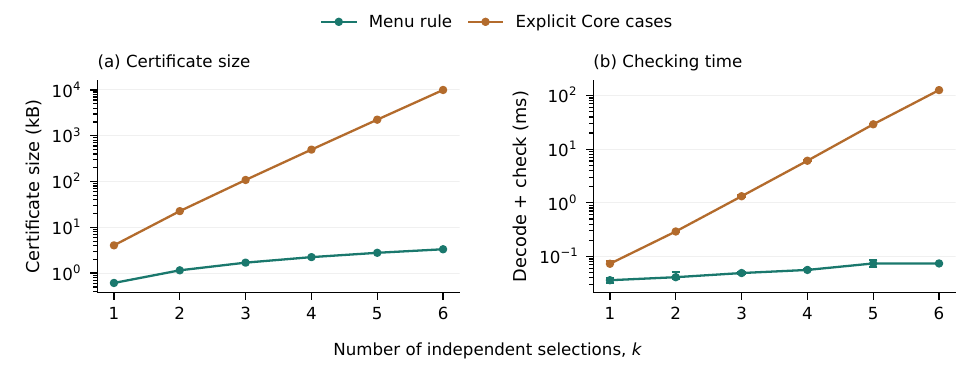}
\caption{\textbf{The same reconstruction encoded with a Menu rule or explicit Core cases.} For each \(k\), the two certificates use identical model pairs, preserve the same input variables, and are checked under identical resource limits. Points show median decoding-and-checking times over five measured runs after one warmup; error bars show the observed minimum and maximum.}
\label{fig:controlled-menu}
\end{figure}

At \(k=6\), the Menu certificate occupies \(3,326\) bytes and takes \(\RPXMenuNativeCheckingMicroseconds\)~\textmu s to decode and check. The explicit Core
cases occupy \(9.99\) MB and take \(\RPXMenuCoreCheckingMilliseconds\)~ms. For this fixed construction, Menu is therefore approximately
\(3,000\) times smaller and \(1,700\) times faster inside the checker.
Figure~\ref{fig:controlled-menu} shows that the gap grows rapidly as independent selections are added.

Both variants encode the same preserved-coordinate reconstruction, isolating the cost of explicit case expansion.
Appendix~\ref{app:controlled-menu} gives the construction and protocol.

\FloatBarrier
\subsection{Validation of the acceptance boundary}

\begin{table}[htbp]
\caption{Validation campaigns and implementation size.}
\label{tab:validation-campaigns}
\centering\small
\renewcommand{\arraystretch}{1.12}
\begin{tabularx}{\textwidth}{@{}P{3.5cm}YY@{}}
\toprule
Validation campaign & Property exercised & Result\\
\midrule
Targeted invalid certificates & Arithmetic rules, case scope, complete output coverage, reconstruction-rule premises, and binary encoding & 12,544 / 12,544 rejected\\[3pt]
Wrong model-pair binding & A certificate must match both source and target models & 240 / 240 rejected\\[3pt]
Deterministic reassembly & RPA compiler reproducibility and stored-artifact consistency & All 240 assembled twice to byte-identical RPX files, matching the stored certificates and passing independent checking\\[3pt]
Exact-import tests & Objective constants, bounds, integrality defaults, ranged rows, malformed inputs & 24 regression tests and 160 independent LP vectors passed\\[3pt]
Implementation size & Size of RPX-owned acceptance implementation & 3,893 production SLOC excluding tests\\
\bottomrule
\end{tabularx}
\end{table}

Table~\ref{tab:validation-campaigns} summarizes the validation campaigns. RPX rejects all \(12,544\) targeted invalid certificates
and all \(240\) tested model-pair mismatches. Import regressions and independent exact-number comparisons test model
interpretation, while repeated assembly tests reproducibility. The RPX-owned acceptance modules contain \(3,893\)
production source lines excluding tests. Appendix~\ref{app:robustness} gives the test construction, results, and counting scope.

\section{Related Work}
\label{sec:related}

\begin{table}[htbp]
\caption{\textbf{Comparison of verification interfaces.} Checked objects, guarantees, and their connection to concrete source--target models.}
\label{tab:interfaces}
\centering\small
\renewcommand{\arraystretch}{1.12}
\begin{tabularx}{\textwidth}{@{}P{2.2cm}P{4.2cm}YP{4.3cm}@{}}
\toprule
Approach & Checked object & Guarantee & Connection to concrete source--target models\\
\midrule
EquivaMap~\cite{zhai2025equivamap} & Formulations, inferred mapping, computed optimum & Mapped optimum passes feasibility and objective tests & Uses extraction and solver-produced solutions\\[3pt]
FLARE~\cite{robbins2026flare} & Lean theorem and proof for parameterized formulations & Universal theorem under formalized semantics & Numerical-model import is a separate bridge\\[3pt]
SOVER~\cite{bhattacharyya2026sover} & Extracted mappings and SMT conditions & Cross-feasibility and objective-order properties & Depends on extracted formulas and SMT encoding\\[3pt]
Certified presolve~\cite{hoen2024presolve} & Reduction trace or proof log & Correctness of logged reductions and recovery & Assumes or records a transformation sequence\\[3pt]
RPX & Two concrete rational MILPs and a certificate & Bidirectional no-worse reconstruction over all feasible solutions & Derives complete checks directly from both models\\
\bottomrule
\end{tabularx}
\end{table}

\subsection{Verification of AI-generated reformulations}

Table~\ref{tab:interfaces} compares the checked objects and guarantees of reformulation-verification systems.
EquivaMap~\cite{zhai2025equivamap} discovers mappings between optimization formulations and, in its implemented verification procedure,
applies the inferred mapping to a computed optimum before checking feasibility and objective agreement. RPX instead
checks a symbolic construction for every feasible solution.

FLARE~\cite{robbins2026flare} uses Lean to prove reformulation theorems for parameterized MILPs. Its guarantees are conditional on
the faithful formalization of the formulations, with import from concrete numerical models forming a separate bridge.
RPX focuses on instance-level certification of concrete source-target pairs. FLARE and RPX also use different
objective relations: FLARE permits strictly increasing objective rescaling, whereas \texttt{MI\_OPT} retains the original
numerical scale and permits pointwise objective improvement in both directions.

SOVER~\cite{bhattacharyya2026sover} separates LLM-based mapping extraction from SMT verification. Its linear path checks cross-feasibility
and a global weak objective-order condition, while its nonlinear path provides a tolerance-aware guarantee. RPX
instead uses exact rational semantics and checks the complete destination-model requirements induced by bidirectional
constructions. These systems therefore occupy complementary points in the design space of reformulation verification.

\subsection{Certification of optimization transformations}

VIPR and rational MIP frameworks certify optimality, bounds, or infeasibility conclusions for a single
model~\cite{cheung2017vipr,eifler2023exact,eifler2024gomory}. Certified MIP presolve records and verifies
reduction steps, including solution recovery for 0-1 models through VeriPB~\cite{hoen2024presolve}. More general dominance and
symmetry proof systems already support solution repair, objective nonincrease, and substantial changes to the
represented search space~\cite{bogaerts2022symmetry,bogaerts2023dominance,vandoornmalen2023proof}.

An RPX certificate establishes the relation between the final rational MILPs directly,
including recovery across different variable representations.

Specialized proof systems for CP, MaxSAT, subgraph algorithms, dynamic programming, and arithmetic
constraints likewise show that domain-specific certificate vocabularies can make complex reasoning
independently checkable~\cite{gocht2020subgraph,gocht2022auditable,berg2023maxsat,demirovic2024states,mcilree2025multiplication}.
Multi-stage proof logging~\cite{flippo2024multistage} and proof-carrying code~\cite{necula1997pcc}
further motivate separating unrestricted artifact production from a restricted acceptance mechanism.

\subsection{Practical and formally trusted checking}

Practical proof logging depends on the logical guarantee, certificate representation, and
checking overhead. VeriPB and CakePB demonstrate high-performance checking with formal guarantees for
pseudo-Boolean optimization~\cite{koops2025practical}, while end-to-end verified subgraph solving combines a specialized checker
with formally verified infrastructure~\cite{gocht2024endtoend}. PBLean uses reflection inside Lean to avoid the memory cost
of constructing large explicit proof terms~\cite{szeider2026pblean}.

\section{Scope and Limitations}
\label{sec:limitations}

\paragraph{Scope of the guarantee.} RPX certifies the relation between the supplied numerical models.
Faithful translation from a natural-language specification requires separate validation.
Exact rational semantics can also yield different conclusions from floating-point solver statuses near
feasibility tolerances. A general executor for applying certified recovery constructions to concrete
solver outputs remains future work.

\paragraph{Trusted implementation.} The acceptance implementation has not been mechanically verified.
The tests in Section~\ref{sec:evaluation} provide implementation evidence; formally verifying the checker or a
translation into a verified backend remains future work.

\paragraph{Coverage of the certificate language.} The Core and fixed reconstruction rules determine the
constructions expressible in RPX. Other correct reformulations may require a different
certificate, additional checking resources, or a new fixed reconstruction rule. Such
extensions retain the same \texttt{MI\_OPT} criterion and complete model checks.

\paragraph{Evaluation scope.} The benchmark measures certification of prepared model pairs and certificates.
End-to-end autonomous reformulation and certificate production remain future work.

\section{Conclusion}
RPX turns the acceptance of an agent-generated MILP reformulation into a fixed, exact model-pair checking
task. Its \texttt{MI\_OPT} criterion requires feasible no-worse constructions in both directions, preserving
feasibility and the optimal objective value while allowing the two models to use different variables
and discard dominated solutions. A small exact Core and fixed reconstruction rules establish the complete
checks generated from both models.

Across \(240\) source-target pairs covering all MIPLIB 2017 Benchmark v2 instances, the same checker and
rule set certify every supplied pair with a median decoding-and-checking time of \(\RPXCheckingMedianSeconds\)~s. These results
show that exact certification of flexible agent-generated reformulations is practical at realistic
MILP scale. More broadly, RPX separates discovery from acceptance: agents may search for reformulations
using arbitrary methods, while adoption of the resulting model rests on independently checkable evidence
tied to the complete source and target models.

\FloatBarrier
\appendix
\counterwithin{table}{section}
\counterwithin{figure}{section}
\numberwithin{equation}{section}
\begin{table}[H]
\caption*{\textbf{Appendix navigation}}
\centering\small
\renewcommand{\arraystretch}{1.1}
\begin{tabularx}{\textwidth}{@{}YP{3.3cm}@{}}
\toprule
Reader question & Location\\
\midrule
What exactly does an RPX certificate mean? & Appendix~\ref{app:directional-semantics}\\
What are the exact Core rules? & Appendix~\ref{app:exact-rules}\\
Why is RPX sound? & Appendix~\ref{app:directional-soundness}\\
What do the six reconstruction rules guarantee? & Appendix~\ref{app:lifts}\\
What does a complete certificate look like? & Appendix~\ref{app:worked-certificates}\\
How were the 240 pairs obtained and measured? & Appendix~\ref{app:evaluation-details}\\
How was the implementation validated? & Appendix~\ref{app:validation-campaigns}\\
How can the results be reproduced? & Appendix~\ref{app:reproduction-workflow}\\
\bottomrule
\end{tabularx}
\end{table}

\section{Formal RPX Semantics and Core Rules}\label{app:core}
\label{app:core-details}

This appendix formalizes the certificate semantics used in Section~\ref{sec:core}.
The fixed reconstruction rules introduced in Section~\ref{sec:lifts} are
specified separately in Appendix~\ref{app:lift-details}.

Throughout, let \(M\) be the input MILP and \(N\) the destination MILP for
one direction of the certificate. All coefficients are exact rationals,
continuous variables range over \(\mathbb{R}\), and declared integer variables
range over \(\mathbb{Z}\). As in Section~\ref{sec:semantics}, both objectives are written as
minimization objectives.

\subsection{Directional certificate semantics}
\label{app:directional-semantics}

A complete RPX certificate contains one directional certificate for
\(M\rightsquigarrow N\) and one for \(N\rightsquigarrow M\). We describe
a single direction.

\paragraph{Core statements and case scopes.} For a normalized affine expression
\(e(x)=a^\top x+b\), the Core records three kinds of statements:
\[
\operatorname{le}(e)
\quad\equiv\quad e\le 0,
\qquad
\operatorname{eq}(e)
\quad\equiv\quad e=0,
\qquad
\operatorname{int}(e)
\quad\equiv\quad e\in\mathbb Z.
\]

A directional certificate starts with the valid statements induced by the
input model \(M\), including its finite variable bounds, finite constraint
sides, and declared integrality conditions. Derived statements may then be
added using the rules in Appendix~\ref{app:exact-rules}.

Case analysis creates nested scopes. A statement derived inside one case is
visible only in that case and its descendants. Sibling cases cannot use one
another's assumptions or derived statements. A reconstruction rule from
Appendix~\ref{app:lift-details} may additionally introduce fresh variables
and fixed relations
among them; these variables and relations are likewise local to the scope
in which they are introduced.

\paragraph{Specifying and checking the destination variables.}
RPX checks a symbolic construction under the assumption that the source
variables \(x\) satisfy \(M\).

For a construction defined by case analysis,
a \emph{terminal case} is a final case that either specifies every destination variable or
is closed by a proved contradiction.
Without case analysis, the whole direction has a single such case.
In each case not closed by contradiction, the certificate gives
an affine expression \(F_j\) for every destination variable \(y_j\),
interpreted as the assignment \(y_j := F_j\).
These expressions may use the source variables and any variables
introduced by fixed reconstruction rules that are available in the
current scope.
Write \(F=(F_1,\ldots,F_n)\) for this vector of formulas, where \(n\)
is the number of destination variables.

RPX substitutes the vector \(F\) for the destination variables
throughout \(N\). The destination model determines all feasibility conditions that must
be established, and \texttt{MI\_OPT} additionally requires
\(f_N(F)\le f_M(x)\).
The certificate then justifies every generated condition.

Suppose \(N\) has variable bounds \(\ell_j\le y_j\le u_j\), row bounds
\(b_r^-\le A_r y\le b_r^+,\)
integer-variable set \(I_N\), and objective \(f_N(y)=c_N^\top y+c_{0,N}\).
For each finite variable bound or constraint side, RPX generates the corresponding check below after substituting \(y:=F\).

\begin{equation*}
\ell_j\le F_j
\quad\Longleftrightarrow\quad
\operatorname{le}(\ell_j-F_j),
\qquad
F_j\le u_j
\quad\Longleftrightarrow\quad
\operatorname{le}(F_j-u_j),
\end{equation*}
\begin{equation*}
b_r^-\le A_rF
\quad\Longleftrightarrow\quad
\operatorname{le}(b_r^- - A_rF),
\qquad
A_rF\le b_r^+
\quad\Longleftrightarrow\quad
\operatorname{le}(A_rF-b_r^+).
\end{equation*}

For every \(j\in I_N\), RPX additionally generates the integrality requirement
\(
\operatorname{int}(F_j)
\).
Finally, RPX generates the no-worse objective requirement
\[
f_N(F)-f_M(x)\le 0,
\]
equivalently,
\[
\operatorname{le}\!\left(
c_N^\top F+c_{0,N}
-
c_M^\top x-c_{0,M}
\right).
\]

Every generated requirement must be justified by a statement visible in the
current scope, subject only to the restricted terminal discharges described
in Appendix~\ref{app:exact-rules}.

A terminal case may instead be closed as infeasible, but only after an exact
contradiction has been derived. A directional certificate is accepted only if
its cases cover every feasible source assignment and each terminal case either
derives an exact contradiction or justifies every generated condition above.

\subsection{Exact Core inference rules}
\label{app:exact-rules}
\label{sec:scalar-rules}

The Core provides a small collection of exact scalar rules. Fixed rules in
Appendix~\ref{app:lift-details} handle more elaborate reconstruction patterns.

\paragraph{Linear and integer inference.}
Let \(p_i\), \(h_j\), and \(z_i\) denote affine expressions whose corresponding statements are visible in the current scope.

\paragraph{Linear inequality.} From visible \(\operatorname{le}(p_i)\) and \(\operatorname{eq}(h_j)\), the Core may derive
\[
\operatorname{le}\left(q_0+\sum_i\alpha_i p_i+\sum_j\beta_j h_j\right),
\]
provided \(q_0\le0\), \(\alpha_i\in\mathbb Q_{\ge0}\), and \(\beta_j\in\mathbb Q\).
Thus valid inequalities may only receive nonnegative multipliers, whereas equalities may receive arbitrary
rational multipliers. The nonpositive rational constant \(q_0\) permits exact weakening by a constant.

\paragraph{Linear equality.} From visible equalities
\(
\operatorname{eq}(h_j)
\),
the Core may derive
\[
\operatorname{eq}\left(\sum_j\beta_j h_j\right),
\qquad\beta_j\in\mathbb Q.
\]
No additional nonzero constant may be introduced.

\paragraph{Integer combination.} From
\(
\operatorname{int}(z_i)
\)
and
\(
\operatorname{eq}(h_j)
\),
the Core may derive
\[
\operatorname{int}\left(k_0+\sum_i k_i z_i+\sum_j\beta_j h_j\right),
\]
where
\(
k_0,k_i\in\mathbb Z,
\beta_j\in\mathbb Q
\).
Rational multiples of equalities are valid because each visible equality establishes \(h_j=0\) in the current scope.

\paragraph{Equality from opposite inequalities.} From
\(
\operatorname{le}(e)
\)
and
\(
\operatorname{le}(-e)
\),
the Core derives
\(
\operatorname{eq}(e)
\).

\paragraph{Integrality-aware rounding.}
Suppose both
\(
\operatorname{le}(a)
\)
and
\(
\operatorname{int}(z)
\)
are visible, and \(a\) and \(z\) have identical variable coefficients. They therefore differ only by a rational constant:
\(
a=z+\delta
\).
Since \(a\le0\) and \(z\) is integer-valued,
\[
z\le-\delta
\quad\Longrightarrow\quad
z\le\lfloor-\delta\rfloor=-\lceil\delta\rceil.
\]
The Core therefore derives
\(
\operatorname{le}\bigl(z+\lceil\delta\rceil\bigr)
\).
This rule rounds only the constant term. It does not round variable coefficients.

\paragraph{Exact weakening from variable bounds.}
Suppose \(\operatorname{le}(p)\) or \(\operatorname{eq}(p)\) is visible and the certificate proposes \(\operatorname{le}(g)\).
The weakening is accepted only when the exact maximum of
\(
g-p
\)
over the currently declared variable box is at most zero. Every nonzero coefficient of \(g-p\) must therefore
have a finite declared bound in the direction required to compute this maximum. If such a bound is unavailable,
the weakening is rejected.
This check uses the declared variable box rather than tighter inequalities derived inside the current case.
Weaken packages a linear-combination argument over these bounds into one instruction:
the checker computes the box maximum directly, avoiding a potentially long list of bound multipliers in the certificate.

\paragraph{Exhaustive case analysis.}
The Core supports exhaustive real and integer case splits. For a rational affine expression \(e\), a real split covers
\(
e\le0
\)
and
\(
e\ge0
\).
For an integer-valued affine expression \(e\), once
\(
\operatorname{int}(e)
\)
has been established, an integer split covers
\(
e\le0
\)
and
\(
e\ge1
\).
There are two forms of case analysis.

\subparagraph{Fact-producing cases.}
The two child cases must derive the same Core statement \(S\). Only \(S\) is exported to the parent scope; all child
assumptions and other derived statements remain local. These internal child derivations use Core scalar rules
and do not invoke reconstruction rules.

\subparagraph{Terminal cases.}
Instead of exporting a single statement, each child may continue independently until it constructs all destination
variables and establishes the complete requirements generated in Appendix~\ref{app:directional-semantics}.
Fixed reconstruction rules from Appendix~\ref{app:lift-details} may be invoked inside such terminal cases.

In both forms, case exhaustiveness is checked by RPX rather than asserted by the certificate.

\paragraph{Contradiction and administrative operations.}
A case may be closed as infeasible only when the current scope contains either
a strictly positive constant inequality asserted to be at most zero, or
a nonzero constant equality asserted to be zero.
These are exact contradictions.
A certificate may explicitly discard a statement that it no longer needs. This operation changes only storage
and visibility; it has no logical effect.

\paragraph{Terminal discharge of trivial requirements.}
Normally, a generated requirement from Appendix~\ref{app:directional-semantics} must match a visible Core statement.
RPX additionally permits a distinguished trivial discharge for the following cases:
\begin{itemize}[leftmargin=*,itemsep=1pt,topsep=3pt]
\item an equality whose expression is exactly the zero constant;
\item an inequality whose expression is a nonpositive constant;
\item an integrality requirement whose expression is an affine integer combination of actual integer variables
currently in scope, with an integer constant.
\end{itemize}
This terminal discharge does not create a reusable Core statement. In particular, it cannot justify a nonconstant inequality.
Table~\ref{tab:logical-core} summarizes the logical Core.
\begin{table}[H]
\caption{The logical Core.}
\label{tab:logical-core}
\centering\small
\renewcommand{\arraystretch}{1.1}
\begin{tabularx}{\textwidth}{@{}P{3.5cm}YP{5.0cm}@{}}
\toprule
Rule & Required information & Result\\
\midrule
Linear inequality & Valid inequalities, equalities, permitted rational multipliers & New valid inequality\\[2pt]
Linear equality & Valid equalities & New valid equality\\[2pt]
Integer combination & Integer-valued expressions and equalities & New integer-valued expression\\[2pt]
Equality pairing & \(e\le0\) and \(-e\le0\) & \(e=0\)\\[2pt]
Round & Bound on an integer-valued affine expression & Rounded constant bound\\[2pt]
Weaken & One valid statement plus finite declared variable bounds & Box-valid weaker inequality\\[2pt]
Real case split & Rational affine expression \(e\) & Exhaustive cases \(e\le0\), \(e\ge0\)\\[2pt]
Integer case split & Integer-valued affine expression \(e\) & Exhaustive cases \(e\le0\), \(e\ge1\)\\[2pt]
Contradiction & Exact inconsistent constant statement & Close infeasible case\\
\bottomrule
\end{tabularx}
\end{table}

\subsection{What the Core can express}
\label{app:core-expressiveness}

Exhaustive case analysis allows different affine constructions in different regions of the source feasible set.

\begin{proposition}[Finite case-based affine encodability]
\label{prop:piecewise}
Let \(T\) be a finite decision tree over \(X_M\). Each internal node performs one of the following exhaustive tests:
(1) for a rational affine expression \(e(x)\),
\(
e(x)\le0
\)
or
\(
e(x)\ge0
\);(2) for an integer-valued affine expression \(e(x)\),
\(
e(x)\le0
\)
or
\(
e(x)\ge1
\).

Each leaf \(\lambda\) assigns a rational affine destination tuple
\(
F_\lambda(x)
\).
Assume that, at every integer split, the required statement
\(
\operatorname{int}(e)
\)
has a finite Core derivation under the assumptions leading to that node. Assume also that, under the assumptions
leading to every feasible leaf \(\lambda\), all destination-model and objective requirements generated for
\(F_\lambda(x)\) have finite Core derivations. Every infeasible leaf must instead have a finite Core derivation
of an exact contradiction.
Then, given sufficient implementation resources, the Core can certify
\(
M\rightsquigarrow N.
\)
\end{proposition}
\begin{proof}
Traverse \(T\) recursively. At each integer node, emit the assumed derivation of \(\operatorname{int}(e)\) and perform
the corresponding exhaustive integer split. At each real node, perform the exhaustive real split. At a feasible
leaf, emit the assumed derivations of all generated destination requirements and return the corresponding tuple
\(F_\lambda(x)\). At an infeasible leaf, emit the assumed contradiction derivation and close the case.
Exhaustiveness ensures that every feasible source solution reaches at least one feasible leaf, and the generated
checks at that leaf establish a feasible destination solution with no larger objective value. Hence the resulting
directional certificate establishes \(M\rightsquigarrow N\).
\end{proof}

The proposition covers, among other patterns, fixed-variable substitution, affine aggregation, bounded integer
enumeration, Boolean case constructions, sign-dependent mappings, and removal of dominated solutions.

\subsection{Directional soundness}
\label{app:directional-soundness}

\begin{theorem}[Directional soundness]
\label{thm:soundness}
Assume exact interpretation of the two models and exact arithmetic. Assume further that the Core rules and
case-scoping semantics of Appendix~\ref{app:exact-rules} are implemented correctly and that every fixed
reconstruction rule used by the certificate correctly instantiates its theorem in Appendix~\ref{app:lift-details}.
If RPX accepts a directional certificate from \(M\) to \(N\), then
\(
M\rightsquigarrow N
\).
\end{theorem}
\begin{proof}
Fix an arbitrary feasible input
\(
x\in X_M
\).
All initial Core statements induced by \(M\) are valid for \(x\). We follow the finite case structure of the
accepted certificate and show inductively that every visible Core statement remains valid for the current
input and for any variables introduced by previously applied reconstruction rules.

The linear inequality and equality rules preserve validity by exact linear combination with the permitted
multiplier signs. Integer combination preserves integer-valuedness. Equality pairing is immediate from the
two opposite inequalities. The rounding rule is valid by integrality of the matched affine expression.
The weakening rule checks the proposed implication exactly over the declared variable box.

At a fact-producing case split, the two children cover the parent domain exhaustively and establish the same
exported statement; therefore that statement is valid in the parent scope. At a terminal split, the children
again cover the parent domain exhaustively, but each continues independently toward a complete destination
construction. Assumptions and derived statements remain within their case scopes.

Whenever a reconstruction rule is invoked, its Core-checked preconditions hold for the current input. By the
corresponding theorem in Appendix~\ref{app:lift-details}, there exists a joint valuation of the introduced
variables satisfying the fixed relations supplied by that rule. Extending the current construction by such a
valuation preserves all previously valid statements and makes the rule's guaranteed relations valid.

A case closed by contradiction contains no feasible realization and therefore requires no destination
construction. Every other terminal case specifies all destination variables. Because RPX generates every
finite destination bound, integrality condition, row side, and the no-worse objective comparison, acceptance
implies that the constructed tuple \(F(x)\) satisfies
\(
F(x)\in X_N
\)
and
\(
f_N(F(x))\le f_M(x)
\).

Thus every feasible \(x\in X_M\) admits a feasible no-worse destination solution. Fixing the certificate's
deterministic case order and one valid choice for each fixed reconstruction rule yields a total map
\(
F:X_M\to X_N
\)
with
\(
f_N(F(x))\le f_M(x)
\text{for all }x\in X_M.
\)
Therefore
\(
M\rightsquigarrow N.
\)
\end{proof}

Applying Theorem~\ref{thm:soundness} to both directional certificates gives
\(
M\rightsquigarrow N
\)
and
\(
N\rightsquigarrow M
\).
By Theorem~\ref{thm:equivalence},
\(
M\simeq_{\mathrm{opt}}N
\).

The soundness argument is independent of which fixed reconstruction rules are available. Adding a new rule
can enlarge the set of certificates that RPX can express compactly, but it does not alter the \texttt{MI\_OPT}
acceptance criterion or remove any of the complete destination-model checks defined in
Appendix~\ref{app:directional-semantics}.

\section{Fixed Reconstruction Rules}
\label{app:lifts}
\label{app:lift-details}

Section~\ref{sec:lifts} introduced six fixed reconstruction rules for cases that are cumbersome to express using only the Core. This appendix gives their exact interfaces and soundness arguments. We use the notation of Appendix~\ref{app:core}:
\(\operatorname{le}(e)\), \(\operatorname{eq}(e)\), and \(\operatorname{int}(e)\) denote, respectively, that an affine expression \(e\) satisfies \(e\le 0\), \(e=0\), or \(e\in\mathbb Z\).
\subsection{Common rule interface}
\label{app:lift-interface}

A reconstruction-rule invocation, called a lift in the RPX instruction format, instantiates a fixed theorem of the form \(A(x)\Longrightarrow\exists z\;R(x,z)\).
Here, \(x\) denotes the variables already available before the invocation. The precondition \(A(x)\) is a fixed collection of affine inequalities, equalities, integrality statements, and structural parameters. The tuple \(z\) contains new variables, and \(R(x,z)\) is the fixed collection of bounds and relations guaranteed for those variables.

The certificate must justify every premise in \(A(x)\) using statements visible in the current case. If all premises are valid, RPX extends that case with the new variables \(z\) and the relations \(R(x,z)\). Subsequent Core steps may use those relations when constructing the destination-model variables.

A reconstruction rule does not remove or modify any existing variable, constraint, or objective coefficient. It only extends the current proof state. When execution leaves a case, variables introduced inside that case are discarded together with its local assumptions.

For batched invocations, every premise is interpreted in the state before the batch. Items in the same batch
may refer to the same existing expressions, but one item may not use variables introduced by another item in
that batch. Variables and relations generated by the completed batch are available to subsequent instructions.

The certificate producer cannot introduce a new theorem or modify the conclusion of an existing rule. Adding a
reconstruction rule requires adding a fixed checker implementation together with a soundness argument.
Regardless of which rules are used, RPX still generates and checks every destination-model bound, integrality condition, constraint, and objective comparison.

\subsection{Local reconstruction rules}
\label{app:local-constructions}

The first three rules introduce variables whose values can be selected locally from established bounds or a finite set of cases. Each admits a direct constructive proof.
\subsubsection{Interval}
\label{app:interval}

Let \(L_1(x),\ldots,L_p(x)\) be affine lower bounds and \(U_1(x),\ldots,U_q(x)\) be affine upper bounds. Either list may be empty.

The Interval rule requires a proof of
\begin{equation}
L_i(x)\le U_j(x)
\qquad
(i=1,\ldots,p,\;j=1,\ldots,q).
\label{eq:interval-premises}
\end{equation}
When one list is empty, no cross-list premise is required. The rule introduces one continuous variable \(z\) and guarantees
\begin{equation}
L_i(x)\le z\quad(i=1,\ldots,p),
\qquad
z\le U_j(x)\quad(j=1,\ldots,q).
\label{eq:interval-bounds}
\end{equation}
\begin{reconstructionlemma}[Interval soundness]
\label{lem:interval}
For every fixed value of \(x\) satisfying the premises in \eqref{eq:interval-premises}, there exists a real number \(z\) satisfying \eqref{eq:interval-bounds}.

\end{reconstructionlemma}
\begin{proof}
If both lists are nonempty, choose \(z=\max_{1\le i\le p}L_i(x)\). This choice satisfies every lower bound. Let \(i^\star\) attain the maximum. For every \(j\), premise \eqref{eq:interval-premises} gives \(z=L_{i^\star}(x)\le U_j(x)\), so every upper bound also holds.

If only the lower-bound list is nonempty, choose its maximum. If only the upper-bound list is nonempty, choose its minimum. If both lists are empty, choose \(z=0\).
\end{proof}

Lemma~\ref{lem:interval} justifies introducing a variable satisfying the generated inequalities once all pairwise premises have been established.
\subsubsection{Menu}
\label{app:menu}

Fix \(K\in\mathbb Z_{\ge 0}\). Let \(A(x)\) and \(I(x)\) be affine expressions available before the rule invocation.

Menu requires
\begin{equation}
A(x)\in\mathbb Z,
\quad
0\le A(x)\le 1,
\quad
I(x)\in\mathbb Z,
\quad
0\le I(x)\le K A(x).
\label{eq:menu-premises}
\end{equation}
The rule introduces variables \(q_0,\ldots,q_K\) and guarantees
\begin{equation}
q_j\in\{0,1\}\quad(j=0,\ldots,K),
\qquad
\sum_{j=0}^{K}q_j=A(x),
\qquad
\sum_{j=0}^{K}j q_j=I(x).
\label{eq:menu-conclusions}
\end{equation}
\begin{reconstructionlemma}[Menu soundness]
\label{lem:menu}
For every fixed value of \(x\) satisfying \eqref{eq:menu-premises}, variables \(q_0,\ldots,q_K\) satisfying \eqref{eq:menu-conclusions} exist.

\end{reconstructionlemma}
\begin{proof}
Because \(A(x)\) is integer-valued and lies in \([0,1]\), \(A(x)\in\{0,1\}\).
If \(A(x)=0\), then \eqref{eq:menu-premises} implies \(I(x)=0\). Set \(q_j=0\) for every \(j\). These values satisfy \eqref{eq:menu-conclusions}.

If \(A(x)=1\), then \(I(x)\) is an integer in \(\{0,\ldots,K\}\). Set \(q_{I(x)}=1\) and set every other \(q_j\) to zero. Then the unweighted sum in \eqref{eq:menu-conclusions} is one, and the weighted sum is \(I(x)\). This argument also covers \(K=0\).
\end{proof}

\subsubsection{Signed}
\label{app:signed}

Fix a rational capacity \(C\ge 0\). Let \(Z(x)\) and \(F(x)\) be affine expressions available before the rule invocation.

Signed requires
\begin{equation}
Z(x)\in\mathbb Z,
\qquad
0\le Z(x)\le 1,
\qquad
-CZ(x)\le F(x)\le CZ(x).
\label{eq:signed-premises}
\end{equation}
The rule introduces one variable \(u\) and guarantees
\begin{equation}
u\in\{0,1\},
\qquad
u\le Z(x),
\qquad
F(x)\le Cu,
\qquad
-F(x)\le C\bigl(Z(x)-u\bigr).
\label{eq:signed-conclusions}
\end{equation}
\begin{reconstructionlemma}[Signed soundness]
\label{lem:signed}
For every fixed value of \(x\) satisfying \eqref{eq:signed-premises}, a variable \(u\) satisfying \eqref{eq:signed-conclusions} exists.

\end{reconstructionlemma}
\begin{proof}
As in Lemma~\ref{lem:menu}, \eqref{eq:signed-premises} implies \(Z(x)\in\{0,1\}\).

If \(Z(x)=0\), then \eqref{eq:signed-premises} forces \(F(x)=0\). Choose \(u=0\). All conclusions follow.

Suppose \(Z(x)=1\). If \(F(x)\ge 0\), choose \(u=1\). Then \(u\le Z(x)\), \(F(x)\le Cu\) follows from \(F(x)\le C\), and the last inequality in \eqref{eq:signed-conclusions} reduces to \(-F(x)\le 0\).

If \(F(x)<0\), choose \(u=0\). Then \(u\le Z(x)\) and \(F(x)\le Cu\) hold, while the last inequality in \eqref{eq:signed-conclusions} follows from \(-F(x)\le C\). The construction also covers \(C=0\) and \(F(x)=0\).
\end{proof}

The rule may be combined with Interval to derive the familiar positive and negative decomposition. Define \(v=Z-u\). The generated relations imply \(u,v\in\{0,1\}\) and \(u+v=Z\). There exist nonnegative continuous variables \(p,n\) satisfying \(p-n=F\), \(p\le Cu\), and \(n\le Cv\).
For example, if \(u=1\), choose \(p=F\) and \(n=0\); if \(v=1\), choose \(p=0\) and \(n=-F\); and if \(Z=0\), choose both variables to be zero. The authoring interface implements this construction by composing Signed with Interval.
\subsection{Network reconstruction rules}
\label{app:network-constructions}

The remaining rules introduce several edge variables that must satisfy shared network constraints. Their soundness follows from standard flow-decomposition and network-integrality arguments.
\subsubsection{Notation}
\label{app:network-notation}

Let \(G=(V,E)\) be a finite directed multigraph. Parallel edges and loops are allowed. For an edge vector \(h\in\mathbb R^E\), define its outgoing-minus-incoming balance at node \(v\) by
\begin{equation*}
(Dh)_v
=
\sum_{\substack{e\in E\\ e\text{ leaves }v}}h_e
-
\sum_{\substack{e\in E\\ e\text{ enters }v}}h_e.
\end{equation*}
A loop appears once in each sum and therefore contributes zero balance.

Fix a source node \(s\). Commodity \(i\) has sink \(t_i\) and nonnegative rational demand \(d_i\). Its required balance vector is \(b^i=d_i(e_s-e_{t_i})\), where \(e_v\) denotes the unit vector of node \(v\). If \(s=t_i\), then \(b^i=0\).
\subsubsection{Flow}
\label{app:flow}

Flow reconstructs destination-specific shipments from an aggregate edge flow while preserving every aggregate edge total.

\begin{reconstructionlemma}[Common-source flow decomposition]
\label{lem:flow}
Let \(G=(V,E)\) be a finite directed multigraph with source \(s\). Let \(k\ge 1\), and for each \(i=1,\ldots,k\), let \(t_i\in V\) and \(d_i\ge 0\). Suppose \(f\in\mathbb R_{\ge 0}^E\) satisfies \(Df=\sum_{i=1}^{k}b^i\). Then there exist nonnegative edge vectors \(g^1,\ldots,g^k\in\mathbb R_{\ge 0}^E\) such that \(Dg^i=b^i\) for each \(i=1,\ldots,k\), and \(\sum_{i=1}^{k}g^i_e=f_e\) for every \(e\in E\).
Repeated sinks, parallel edges, loops, disconnected components, and zero demands are allowed.

\end{reconstructionlemma}
\begin{proof}
Maintain a nonnegative residual flow \(r\), initially \(r=f\), together with the unsatisfied demand of every commodity. Commodities with \(s=t_i\) have zero balance and require no path flow.

Consider a commodity with positive unsatisfied demand. We first show that its sink is reachable from \(s\) through edges with positive residual flow. Let \(S\subseteq V\) be the set of nodes reachable from \(s\) through positive-residual edges. If the sink were outside \(S\), no positive-residual edge would leave \(S\). Therefore,
\begin{equation}
\sum_{v\in S}(Dr)_v\le 0.
\label{eq:flow-cut-balance}
\end{equation}
On the other hand, the residual balance over \(S\) equals the total unsatisfied demand of commodities whose sinks lie outside \(S\), which is positive. This contradicts \eqref{eq:flow-cut-balance}. Hence a positive-residual \(s\)-to-sink path exists.

Remove cycles from the path so that it is simple. Send along it the smaller of the path bottleneck and the commodity's remaining demand. Subtract this amount from the residual and add it to the commodity flow. Each iteration either satisfies one remaining demand or reduces at least one positive residual edge to zero. A zero residual edge never becomes positive again, so finitely many iterations satisfy every demand.

After all demands have been met, the remaining residual \(r\) has zero balance, \(Dr=0\). Because \(k\ge 1\), add this residual circulation to any one commodity flow. This preserves that commodity's balance and makes the edgewise sum of the commodity flows equal to \(f\). Nonnegativity is maintained throughout.
\end{proof}

\paragraph{RPX instantiation.}
Flow takes a directed multigraph \(G=(V,E)\), a source \(s\), sinks \(t_1,\ldots,t_k\) with exact nonnegative rational demands \(d_1,\ldots,d_k\), and an old-scope affine expression \(f_e(x)\) for each edge. The certificate proves \(f_e(x)\ge 0\) on every edge and \((Df(x))_v=\sum_{i=1}^{k}b^i_v\) at every node. RPX introduces nonnegative continuous variables \(g^i_e\) satisfying \(Dg^i=b^i\) for each commodity and \(\sum_{i=1}^{k}g^i_e=f_e(x)\) on every edge. Any objective comparison or additional side constraints on these flows must be established through the destination-model checks.
\subsubsection{Subflow}
\label{app:subflow}

Subflow removes excess circulation from an aggregate flow while preserving all required node balances.
Let \(B=\sum_{i=1}^{k}d_i\) denote the total demand.

\begin{reconstructionlemma}[Demand-bounded subflow]
\label{lem:subflow}
Let \(f\in\mathbb R_{\ge 0}^E\) satisfy \(Df=b\), where \(b=\sum_{i=1}^{k}b^i\). Then there exists \(g\in\mathbb R^E\) satisfying \(Dg=b\) and, for every edge \(e\), \(0\le g_e\le f_e\) and \(g_e\le B\).
Loops, parallel edges, repeated sinks, disconnected components, and zero demands are allowed.

\end{reconstructionlemma}
\begin{proof}
Apply the path-extraction procedure from the proof of Lemma~\ref{lem:flow}, but discard the final zero-balance residual instead of assigning it to a commodity. Let \(g\) be the sum of all extracted path flows.

Each extracted path carries flow that was available in the current residual, so the total extracted amount on any edge does not exceed \(f_e\). Hence \(0\le g_e\le f_e\). The extracted paths meet every nonzero source-to-sink demand, so their sum has balance \(b\).

Each path can be chosen simple. It therefore uses any fixed edge at most once. The total weight of all extracted paths is at most the total demand \(B\), implying \(g_e\le B\) for every edge. A demand whose sink is the source contributes zero balance and requires no path. A loop cannot appear in a simple extracted path.
\end{proof}

\paragraph{RPX instantiation.}
Subflow takes the same input data as Flow, with \(t_i\ne s\). The certificate proves \(f_e(x)\ge 0\) on every edge and \(Df(x)=b\). RPX introduces continuous variables \(g_e\) satisfying \(0\le g_e\le B\), \(Dg=b\), and \(g_e\le f_e(x)\) on every edge. The rule may remove zero-balance circulation and gives no objective guarantee; any required objective comparison remains part of the complete destination-model checks.

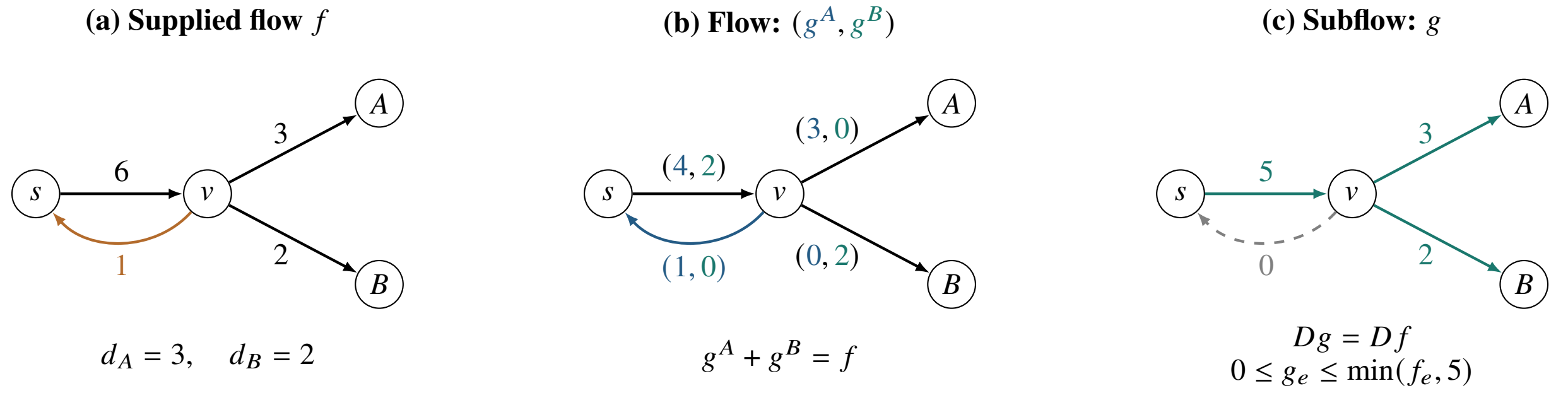
\begin{figure}[htbp]
\centering
\begin{tikzpicture}[font=\small,
vertex/.style={circle,draw,fill=white,inner sep=1.5pt,minimum size=5mm},
edge/.style={-{Latex[length=1.5mm]},thick}]
\foreach \shift/\prefix in {0/f,6.0/d,12.0/g}{
  \node[vertex] (\prefix s) at (\shift,0) {$s$};
  \node[vertex] (\prefix v) at (\shift+1.8,0) {$v$};
  \node[vertex] (\prefix a) at (\shift+3.6,0.95) {$A$};
  \node[vertex] (\prefix b) at (\shift+3.6,-0.95) {$B$};
}
\node[font=\small\bfseries] at (1.8,1.8) {(a) Supplied flow $f$};
\draw[edge] (fs)--node[above] {6}(fv);
\draw[edge,rpxorange] (fv) to[bend left=48] node[below] {1} (fs);
\draw[edge] (fv)--node[above left,inner sep=1pt] {3}(fa);
\draw[edge] (fv)--node[below left,inner sep=1pt] {2}(fb);
\node[align=center] at (1.8,-1.7) {$d_A=3,\quad d_B=2$};

\node[font=\small\bfseries] at (7.8,1.8) {(b) Flow: $(\textcolor{rpxblue}{g^A},\textcolor{rpxteal}{g^B})$};
\draw[edge] (ds)--node[above] {$(\textcolor{rpxblue}{4},\textcolor{rpxteal}{2})$}(dv);
\draw[edge,rpxblue] (dv) to[bend left=48] node[below] {$(\textcolor{rpxblue}{1},\textcolor{rpxteal}{0})$} (ds);
\draw[edge] (dv)--node[above left,inner sep=1pt] {$(\textcolor{rpxblue}{3},\textcolor{rpxteal}{0})$}(da);
\draw[edge] (dv)--node[below left,inner sep=1pt] {$(\textcolor{rpxblue}{0},\textcolor{rpxteal}{2})$}(db);
\node[align=center] at (7.8,-1.7) {$g^A+g^B=f$};

\node[font=\small\bfseries] at (13.8,1.8) {(c) Subflow: $g$};
\draw[edge,rpxteal] (gs)--node[above] {5}(gv);
\draw[edge,gray,dashed] (gv) to[bend left=48] node[below] {0} (gs);
\draw[edge,rpxteal] (gv)--node[above left,inner sep=1pt] {3}(ga);
\draw[edge,rpxteal] (gv)--node[below left,inner sep=1pt] {2}(gb);
\node[align=center] at (13.8,-1.7) {$Dg=Df$\\$0\le g_e\le\min(f_e,5)$};
\end{tikzpicture}
\caption{\textbf{Flow and Subflow solve different reconstruction problems.} Flow decomposes an aggregate flow into destination-specific flows while preserving the total on every edge. In (b), commodity \(A\) carries the unit circulation \(s\to v\to s\), giving \((4,2)\) on \(s\to v\) and \((1,0)\) on \(v\to s\) while delivering three units to \(A\). Subflow may remove circulation, but preserves the required node balances and never exceeds the supplied edge flow or total demand.}
\label{fig:flow-subflow}
\end{figure}

Figure~\ref{fig:flow-subflow} illustrates the distinction. Flow preserves every edge total but separates it by destination; Subflow may reduce edge values by discarding circulation while preserving the required deliveries.

\FloatBarrier

\subsubsection{UnitFlowRound}
\label{app:unitflow}

UnitFlowRound reconstructs binary edge decisions from fractional edge values while preserving all node balances and not increasing a fixed linear cost.

\begin{reconstructionlemma}[Cost-nonincreasing binary unit-flow existence]
\label{lem:unit-flow-round}
Let \(G=(V,E)\) be a finite directed multigraph with \(|V|>0\). Let \(D\) be its outgoing-minus-incoming incidence operator, and let \(c\in\mathbb Q^E\) be fixed. For every \(f\in[0,1]^E\) satisfying \(Df\in\mathbb Z^V\), there exists \(g\in\{0,1\}^E\) such that \(Dg=Df\) and \(c^\top(g-f)\le 0\).
Loops, parallel edges, disconnected graphs, isolated nodes, negative costs, and \(E=\varnothing\) are allowed.

\end{reconstructionlemma}
\begin{proof}
Call an edge fixed if its current value is \(0\) or \(1\). Initially freeze every fixed edge and consider the undirected multigraph induced by the remaining strictly fractional edges.

Suppose this multigraph contains a cycle. Choose a signed unit circulation \(d\) supported on that cycle. For an edge aligned with one orientation of the cycle, set its component in \(d\) to \(1\); for an oppositely oriented edge, set it to \(-1\). A loop gives a one-edge cycle, and parallel edges may give a two-edge cycle. In every case, \(Dd=0\). Replace \(d\) by \(-d\) if necessary so that \(c^\top d\le 0\). Because every edge in the cycle is strictly inside \([0,1]\), there exists \(\alpha>0\) such that \(f+\alpha d\in[0,1]^E\) and at least one additional edge reaches \(0\) or \(1\). Choose the largest such \(\alpha\). This update preserves the node balances because \(Dd=0\), and it does not increase the cost because
\[
c^\top(f+\alpha d)-c^\top f
=
\alpha c^\top d
\le 0.
\]
Freeze every newly integral edge and repeat. Each iteration fixes at least one additional edge, so the procedure terminates after at most \(|E|\) iterations.

At termination, the undirected graph induced by the remaining fractional edges is a forest. Suppose this forest is nonempty and take a leaf node. Contributions from all fixed edges are integer because their values are binary. The total node balance is also integer by assumption. The leaf is incident to exactly one remaining fractional edge, so that edge's signed contribution must be integer. This contradicts its value being strictly between \(0\) and \(1\).

Therefore no fractional edge remains. The resulting vector \(g\) is binary, has the same node balances as \(f\), and has no greater cost. The argument applies independently to disconnected components; isolated nodes and an empty edge set require no update.
\end{proof}

This is the unit-capacity form of the standard integrality property of minimum-cost network flow~\cite{ahuja1993network}.

\paragraph{RPX instantiation.}
UnitFlowRound takes a directed multigraph with \(n>0\) nodes and \(m\ge 0\) ordered edges, an old-scope affine expression \(f_e(x)\) for each edge, and a fixed exact rational cost vector \(c\). The certificate proves \(0\le f_e(x)\le 1\) on every edge and \((Df(x))_v\in\mathbb Z\) at every node. RPX introduces integer variables \(g_e\) with \(0\le g_e\le 1\), hence binary, and generates
\begin{align}
Dg&=Df(x),
\label{eq:unitflow-generated-balance}\\
\sum_{e\in E}c_e\bigl(g_e-f_e(x)\bigr)&\le 0.
\label{eq:unitflow-generated-cost}
\end{align}
The cost vector is stored sparsely as a strictly increasing list of distinct edge indices with nonzero exact rational coefficients; omitted edges have zero cost.

The generated relations preserve node balances and bound the cost, as specified in
\eqref{eq:unitflow-generated-balance}--\eqref{eq:unitflow-generated-cost}.
\subsection{Soundness and integration with the Core}
\label{app:reconstruction-soundness}

Lemmas~\ref{lem:interval}-\ref{lem:unit-flow-round} establish the extension property required by Appendix~\ref{app:directional-soundness}:
whenever a rule's premises hold, its introduced variables can satisfy
all generated relations. During checking, RPX verifies the premises and introduces the guaranteed relations into the proof state.

\section{Worked Certificates}
\label{app:worked-certificates}

Appendix~\ref{app:core} defines the exact semantics of directional certificates and the Core rules. Appendix~\ref{app:lifts} specifies the six fixed reconstruction rules. This appendix combines these components in two worked examples.

Section~\ref{app:running-certificate} gives an executable certificate for the same two-link running example used in Sections~\ref{sec:semantics}--\ref{sec:method}. It illustrates exhaustive case analysis, integer rounding, and complete destination-model checks using only the Core. Section~\ref{app:waihi} presents a benchmark reformulation in which UnitFlowRound reconstructs hundreds of binary choices jointly.

The companion artifact contains the concrete models, executable RPA sources, binary RPX certificates, and canonical disassemblies. Appendix~\ref{app:reproduction-workflow} gives the reproduction workflow.

\subsection{Running two-link example}
\label{app:running-certificate}
\label{sec:assembly-example}

We use exactly the source and target models from Section~\ref{sec:example}:
\begin{equation}
\begin{aligned}
M:\quad \min\ &a+2b-2q\\
\text{s.t.}\quad &q\le2a,\quad q\le2b,\\
&q\in\{0,1,2\},\quad a,b\in\{0,1\},
\end{aligned}
\label{eq:worked-running-source}
\end{equation}
and
\begin{equation}
\begin{aligned}
N:\quad \min\ &3z-2q\\
\text{s.t.}\quad &q\le2z,\\
&q\in\{0,1,2\},\quad z\in\{0,1\}.
\end{aligned}
\label{eq:worked-running-target}
\end{equation}
The source has two opening decisions \(a,b\); the target contracts them into one decision \(z\). Both constructions retain \(q\).

\subsubsection{Constructions and return checks}
\label{app:running-forward}
\label{app:running-reverse}

The constructions are those of Sections~\ref{sec:semantics}--\ref{sec:method}. Both retain \(q\): forward chooses \(z=0\) when \(q=0\) and \(z=1\) when \(q\ge1\), satisfying \(3z-a-2b\le0\); reverse sets \(a=b=z\), so \(a+2b-3z=0\). Each of the two forward returns binds eight target-model checks, while the reverse return binds twelve source-model checks. Table~\ref{tab:running-checks} maps the source lines in Listing~\ref{lst:running-rpa} to these checks.

\Needspace{29\baselineskip}
\subsubsection{Executable RPA certificate}
\label{app:running-authoring}

Listing~\ref{lst:running-rpa} is the complete executable RPA certificate, including both directions and every return binding. In the forward block, \texttt{x0,x1,x2} denote \(q,a,b\), and \texttt{r0.hi,r1.hi} denote the two source capacity inequalities. In reverse, coordinates restart at the target: \texttt{x0,x1} denote \(q,z\). Comments identify the arithmetic consequences.

\begin{lstlisting}[
  caption={Complete executable RPA certificate for the running two-link example.},
  label={lst:running-rpa},
  captionpos=b,
  basicstyle=\ttfamily\fontsize{9}{10.4}\selectfont,
  frame=tb,
  breaklines=true,
  columns=fullflexible,
  keepspaces=true,
  morecomment={[l]{;}},
  commentstyle=\color{black!60}\itshape,
  numbers=left,
  numberstyle=\tiny\color{black!60},
  numbersep=7pt,
  xleftmargin=1.8em
]
RPA 13
.forward
  split x0, int v0.int  ; q=0 or q>=1, using integer q
  left a0
    lin.le f0, v1.lo + 2*v2.lo  ; -a-2b <= 0
    ret (x0, 0) (v0.lo, v0.hi, v0.int, exact, exact, exact, a0, f0)
  right a0
    lin.le f1, 1/2*a0 + 1/2*r0.hi  ; 1/2-a <= 0
    lin.int f2, -v1.int
    round f3, f1, f2  ; 1-a <= 0
    lin.le f4, 1/2*a0 + 1/2*r1.hi  ; 1/2-b <= 0
    lin.int f5, -v2.int
    round f6, f4, f5  ; 1-b <= 0
    lin.le f7, f3 + 2*f6  ; 3-a-2b <= 0
    ret (x0, 1) (v0.lo, v0.hi, v0.int, exact, exact, exact, v0.hi, f7)
  end
.reverse
  ret (x0, x1, x1) (
    v0.lo, v0.hi, v0.int, v1.lo, v1.hi, v1.int,
    v1.lo, v1.hi, v1.int, r0.hi, r0.hi, exact)
\end{lstlisting}

\begin{table}[H]
\caption{Source-line bindings for the generated model checks in Listing~\ref{lst:running-rpa}.}
\label{tab:running-checks}
\centering\small
\renewcommand{\arraystretch}{1.1}
\begin{tabularx}{\textwidth}{@{}P{1.6cm}P{5.9cm}Y@{}}
\toprule
Lines & Generated model checks & Evidence\\
\midrule
6, 15 & \(0\le q\le2\), \(q\in\mathbb Z\) & The first three bindings retain the source domain of \(q\).\\
6, 15 & \(z\ge0\), \(z\le1\), \(z\in\mathbb Z\) & The chosen constants \(0\) and \(1\) satisfy all three checks.\\
6, 15 & \(q-2z\le0\) & The left assumption gives \(q\le0\); the right return uses \(q\le2\).\\
5--15 & \(3z-a-2b\le0\) & Lower bounds prove the left objective check; two integer roundings prove the right one.\\
18--20 & \(0\le q\le2\), \(q\in\mathbb Z\) & Copy the target coordinate \(q\).\\
18--20 & \(0\le a,b\le1\), \(a,b\in\mathbb Z\) & Both reconstructed coordinates copy the binary target variable \(z\).\\
18--20 & \(q-2a\le0\), \(q-2b\le0\) & Bind the same target capacity inequality to both source rows.\\
18--20 & \(a+2b-3z\le0\) & Substitution gives zero, so the objective check is an exact equality.\\
\bottomrule
\end{tabularx}
\end{table}

The companion artifact provides both models, the RPA source, the assembled RPX certificate, and its canonical disassembly. RPX accepts the certificate, and reassembling the disassembly reproduces the same bytes.

\subsection{Joint reconstruction in \texttt{neos-4532248-waihi}}
\label{app:waihi}

\subsubsection{Model pair}
\label{app:waihi-models}
The instance contains 145 selection blocks, each with 574 binary choice
variables \(w_j\), binary availability variables \(x_j\), and a binary
activation variable \(y\). The blocks share a nonnegative variable \(t\),
with objective \(\min t\). For one block, the source requires
\begin{equation}
\begin{aligned}
\sum_j w_j &= y, \qquad w_j\in\{0,1\},\quad w_j\le x_j,\\
t-642y-d_jw_j &\ge -642 \qquad \text{for every }j,
\end{aligned}
\label{eq:waihi-source}
\end{equation}
where each \(d_j\) is a positive rational.
The target retains the selection and availability constraints and the
domains of \(x_j,y,t\), relaxes the choices to \(0\le w_j\le x_j\), and
replaces the individual cost rows by
\begin{equation}
t\ge \sum_j d_jw_j.
\label{eq:waihi-target}
\end{equation}
Across all blocks, 83,230 rows are replaced by 145 aggregate rows.

\subsubsection{Forward construction}
\label{app:waihi-forward}

The forward construction retains every variable value. If \(y=0\), all
\(w_j=0\), so \eqref{eq:waihi-target} follows from \(t\ge0\). If \(y=1\),
exactly one \(w_{j^\star}=1\), and its source row gives
\(t\ge d_{j^\star}=\sum_j d_jw_j\). Thus the target is feasible with
unchanged objective.

\subsubsection{Network representation}
\label{app:waihi-network}

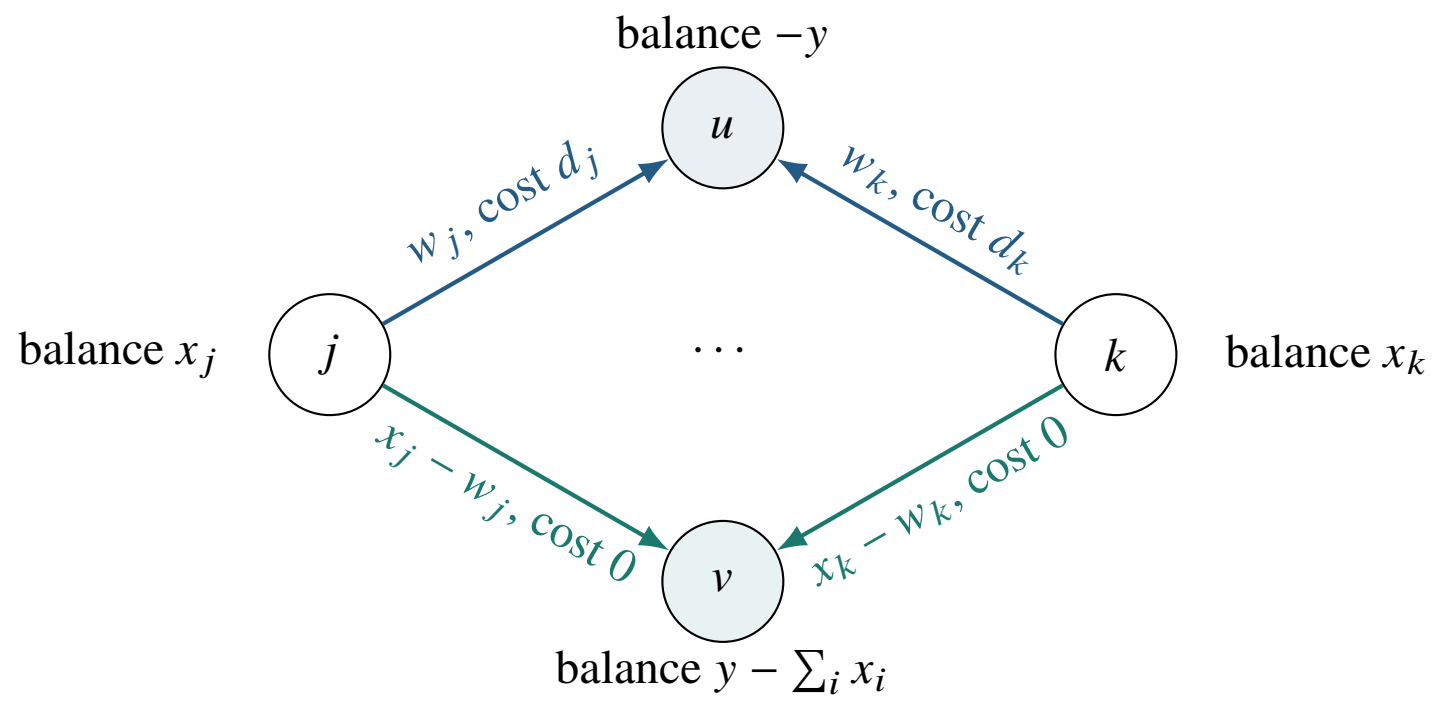
\begin{figure}[htbp]
\centering
\begin{tikzpicture}[
  font=\small,
  arr/.style={-{Latex[length=2mm]},thick},
  vertex/.style={circle,draw,fill=white,minimum size=8mm,inner sep=1pt}]
\node[vertex] (j) at (-2.6,0) {$j$};
\node[vertex] (k) at (2.6,0) {$k$};
\node at (0,0) {$\cdots$};
\node[vertex,fill=rpxblue!10] (u) at (0,1.5) {$u$};
\node[vertex,fill=rpxteal!10] (v) at (0,-1.5) {$v$};
\draw[arr,rpxblue] (j)--node[above,sloped] {$w_j$, cost $d_j$}(u);
\draw[arr,rpxblue] (k)--node[above,sloped] {$w_k$, cost $d_k$}(u);
\draw[arr,rpxteal] (j)--node[below,sloped] {$x_j-w_j$, cost $0$}(v);
\draw[arr,rpxteal] (k)--node[below,sloped] {$x_k-w_k$, cost $0$}(v);
\node[anchor=east] at (-3.2,0) {balance $x_j$};
\node[anchor=west] at (3.2,0) {balance $x_k$};
\node at (0,2.1) {balance $-y$};
\node at (0,-2.1) {balance $y-\sum_i x_i$};
\end{tikzpicture}
\caption{\textbf{Network encoding of one selection block.} For each choice \(j\), edge \(j\to u\) carries \(w_j\) with cost \(d_j\), and edge \(j\to v\) carries \(x_j-w_j\) with zero cost. UnitFlowRound makes the edge values binary while preserving node balances, yielding \(w'_j\le x_j\), \(\sum_j w'_j=y\), and \(\sum_j d_jw'_j\le\sum_j d_jw_j\).}
\label{fig:selection-block}
\label{fig:onehot-network}
\end{figure}

The reverse direction must recover binary choices while preserving
selection, availability, and the weighted-cost bound.
Figure~\ref{fig:selection-block} shows the network encoding.
For each block, construct a network with choice nodes \(j\) and two
additional nodes \(u,v\). Edge \(j\to u\) carries \(w_j\) with cost \(d_j\);
edge \(j\to v\) carries \(x_j-w_j\) with zero cost.
All edge values lie in \([0,1]\). Their outgoing-minus-incoming balances
are \(x_j\) at node \(j\), \(-y\) at \(u\), and
\(y-\sum_j x_j\) at \(v\), hence integer-valued.
This 576-node, 1,148-edge network therefore satisfies the premises of
UnitFlowRound.

\FloatBarrier

The rule supplies binary edge variables \(w'_j\) on \(j\to u\) and
\(h'_j\) on \(j\to v\), with
\begin{equation}
\begin{gathered}
w'_j,h'_j\in\{0,1\},\qquad w'_j+h'_j=x_j,\\
\sum_j w'_j=y,\qquad
\sum_j d_jw'_j\le\sum_j d_jw_j\le t.
\end{gathered}
\label{eq:waihi-recovery}
\end{equation}
These relations recover the selection and availability constraints,
including \(w'_j\le x_j\), without increasing weighted cost.

\subsubsection{Complete source-model checks}
\label{app:waihi-source-recovery}
Replace \(w_j\) by \(w'_j\), retaining \(x_j,y,t\) and all variables
outside the selection blocks.
Since \(d_j>0\), \(w'_j\ge0\), and \(y\le1\), the cost bound in
\eqref{eq:waihi-recovery} implies, for every \(j\),
\begin{equation}
t-642y-d_jw'_j
\ge t-642y-\sum_k d_kw'_k
\ge -642y
\ge -642.
\label{eq:waihi-source-recovery}
\end{equation}
Thus every original cost row holds, and the objective remains \(t\).
The certificate applies this construction to all 145 blocks.
RPX then generates and checks every bound, integrality condition,
constraint, and objective comparison for the complete recovered source solution.

\section{Evaluation Details}
\label{appD:evaluation}
\label{app:evaluation-details}
\label{app:coverage}
This appendix describes the evaluated model pairs, measurement protocol,
supplementary statistics, and controlled Menu experiment. Implementation
validation and artifact reproduction are reported in Appendix~\ref{app:robustness}.

\subsection{Evaluation set and scope}
\label{appD:scope}
The evaluation comprises 240 source--target pairs, one for each MIPLIB 2017
Benchmark v2 instance. An AI agent produced the target formulations under a
30-minute execution budget, with instructions to justify its changes or report
that no justified improvement was found. Certificates were prepared after target selection.

We evaluate the coverage and checking performance of RPX on these pairs using
one common Core and the same six reconstruction rules.

\subsection{Measurement protocol}
\label{appD:protocol}
\label{app:protocol}
Coverage and timing use the same fixed configuration: $10^{10}$ charged work units,
8,000,000 simultaneously stored statements, and case depth 64.
All 240 pairs are accepted. A charged work unit counts a deterministic implementation
operation, such as an expression-term traversal or a statement-table probe.

Timing uses an Apple M5 with 16\,GiB of memory. Cases run sequentially,
with one excluded warmup followed by five fresh-process measurements;
filesystem caches are warm and scheduling is unpinned. All five measured
runs of every pair return Accepted. Reported times are per-pair medians;
peak RSS is the largest observed process peak for each pair.

Inputs are prepared canonical exact models and binary certificates. Agent
execution, certificate production, text-model import, downstream solving,
and numerical solution recovery are excluded. Table~\ref{tab:appD:timing}
defines the timing boundaries. Decode + check is the primary metric used
in Section~\ref{sec:evaluation}.

\begin{table}[htbp]
\centering
\small
\caption{Timing boundaries. Totals denote the sum of per-pair medians across all 240 pairs.}
\label{tab:appD:timing}
\begin{tabular}{@{}lp{0.60\linewidth}r@{}}
\toprule
Metric & Included work & Total \\
\midrule
Kernel & Certificate execution and generated model checks & \RPXKernelTotalSeconds\,s \\
Decode + check & Binary certificate decoding and kernel execution & \RPXCheckingTotalSeconds\,s \\
CLI total & Complete measured work inside the RPX command & \RPXCLITotalTotalSeconds\,s \\
External wall & Elapsed time around the fresh RPX process & \RPXWallTotalSeconds\,s \\
\bottomrule
\end{tabular}
\end{table}

\subsection{Supplementary model and certificate statistics}
\label{appD:statistics}
Table~\ref{tab:appD:rules} supplements the mechanism counts in Figure~\ref{fig:mechanism-usage} with
the number of variables introduced by reconstruction rules, summed over all calls and cases.
Certificate counts overlap. Of 149 UnitFlowRound invocations, 145 occur in
\texttt{neos-4532248-waihi}; the remaining four occur in \texttt{neos-960392},
\texttt{rail02}, \texttt{supportcase6}, and \texttt{tbfp-network}.

\begin{table}[htbp]
\centering
\small
\caption{Reconstruction-rule use. Introduced variables are summed over all calls and cases.}
\label{tab:appD:rules}
\begin{tabular}{@{}lrr@{}}
\toprule
Rule & Certificates & Introduced variables \\
\midrule
Interval & 42 & 424,216 \\
Menu & 20 & 1,083,853 \\
Signed & 22 & 343,352 \\
Flow & 1 & 2,350 \\
Subflow & 2 & 50,726 \\
UnitFlowRound & 5 & 640,430 \\
\bottomrule
\end{tabular}
\end{table}

\begin{figure}[htbp]
\centering
\begin{minipage}[t]{0.49\linewidth}
\centering
\includegraphics[width=\linewidth]{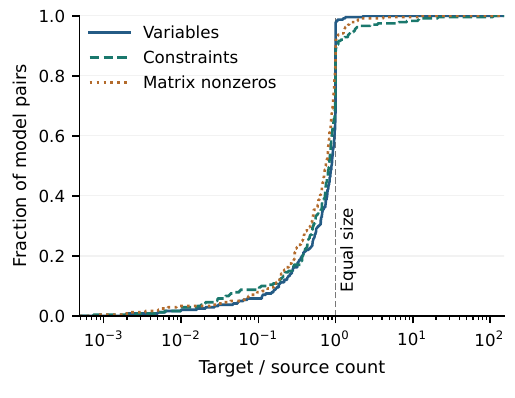}
\par\smallskip
{\small (a) Changes in model size}
\end{minipage}\hfill
\begin{minipage}[t]{0.49\linewidth}
\centering
\includegraphics[width=\linewidth]{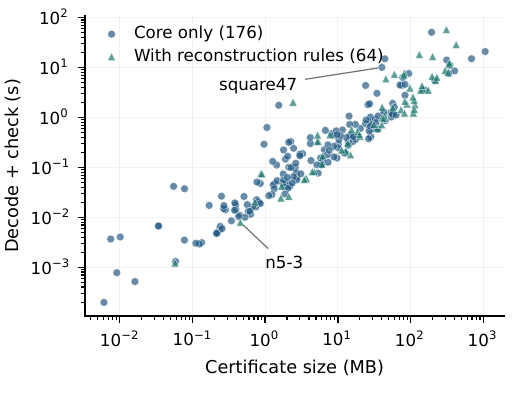}
\par\smallskip
{\small (b) Certificate size and checking time}
\end{minipage}
\caption{Supplementary statistics for all 240 model pairs.
(a) Empirical cumulative distributions of target/source variable, constraint,
and matrix-nonzero counts; the vertical line at one denotes unchanged size.
(b) Certificate size and median decoding-and-checking time for each pair,
grouped by use of reconstruction rules.}
\label{fig:appD:diagnostics}
\end{figure}

Figure~\ref{fig:appD:diagnostics}(a) shows the distributions of target/source
model-size ratios. Targets have more variables, constraints, and matrix
nonzeros than their sources in 6, 32, and 21 pairs, respectively.
Figure~\ref{fig:appD:diagnostics}(b) relates certificate size to checking time
and distinguishes Core-only certificates from those using reconstruction
rules. These distributions describe variation across the benchmark pairs.
Section~\ref{sec:checking-cost} and Table~\ref{tab:verification240} report aggregate cost quantiles.

\subsection{Controlled Menu experiment}
\label{appD:menu}
\label{app:controlled-menu}
\paragraph{Model family and fixed construction.}
For $k=1,\ldots,6$ independent selections, each instantiating Menu with $K=2$, the compact model has
\[
  A_i\in\{0,1\},\qquad I_i\in\mathbb Z,\qquad
  0\le I_i\le 2A_i,\qquad i=1,\ldots,k.
\]
The expanded model retains these variables and adds binary
$q_{i0},q_{i1},q_{i2}$ satisfying
\[
  q_{i0}+q_{i1}+q_{i2}=A_i,\qquad q_{i1}+2q_{i2}=I_i.
\]
Both models minimize $-\sum_i I_i$. The compact-to-expanded construction
preserves $(A,I)$ and supplies their uniquely determined option variables;
the reverse construction projects away $q$.

\paragraph{Certificate variants and protocol.}
The Menu certificate uses one batched invocation. The explicit Core-case certificate
splits each $(A_i,I_i)$ into its four feasible states,
$(0,0),(1,0),(1,1),(1,2)$, producing $4^k$ combinations and repeating
subsequent checks in each case. Both variants therefore certify the same
model pair and the same construction using the same production checker.

All six explicit case expansions complete under construction ceilings of
256\,MiB, one million generated instructions, depth 60, and 60\,s.
Both variants use identical checker limits. Timing follows
Appendix~\ref{appD:protocol}, with Menu-rule and explicit-case runs alternating in order.

\paragraph{Results and interpretation.}
Figure~\ref{fig:controlled-menu} reports all six values of $k$.
Table~\ref{tab:appD:menu} adds timing ranges and fresh-process times at $k=6$.
The Menu certificate is approximately 3,000 times smaller and 1,700 times
faster to decode and check than this explicit case expansion.

\begin{table}[htbp]
\centering
\small
\caption{Controlled comparison at $k=6$. Brackets give the observed
minimum--maximum over five measured runs; MB is decimal.}
\label{tab:appD:menu}
\begin{tabular}{@{}llrr@{}}
\toprule
Representation & Size & Decode + check (ms) & Process time (ms) \\
 & & Median [range] & Median \\
\midrule
Menu rule & 3,326 bytes & \RPXMenuNativeCheckingMilliseconds{} [\RPXMenuNativeMinMilliseconds--\RPXMenuNativeMaxMilliseconds] & \RPXMenuNativeProcessMilliseconds \\
Explicit Core cases & 9.99 MB & \RPXMenuCoreCheckingMilliseconds{} [\RPXMenuCoreMinMilliseconds--\RPXMenuCoreMaxMilliseconds] & \RPXMenuCoreProcessMilliseconds \\
\bottomrule
\end{tabular}
\end{table}

Generated models, certificates, and raw measurements accompany
the reproduction artifact described in Appendix~\ref{app:reproduction-workflow}.

\section{Implementation Validation and Reproducibility}
\label{appE:implementation}
\label{app:robustness}
This appendix describes input interpretation, implementation tests, code size,
and the reproduction workflow.

\subsection{Accepted input semantics}
\label{appE:semantics}
\label{app:input-semantics}
\label{app:import-contract}
RPX interprets integers, finite decimals, scientific notation, and integer
fractions as exact rationals; continuous variables range over the reals,
and declared integer and binary domains are enforced without feasibility
tolerances. Import resolves the default bounds and declarations of the
supported MPS/LP formats; an absent bound side in the resulting canonical
model denotes infinity. Duplicate coefficients are summed exactly, and
maximization is normalized by negating the complete affine objective,
including its constant. Contradictory finite bounds remain a description
of an infeasible model.

When checking starts from MPS/LP text, import is part of the acceptance boundary.
The importer supports linear and ranged rows, explicit bounds, and integrality
declarations. It rejects ambiguous syntax, nonlinear terms, SOS and indicator
constraints, and multiple RHS, range, or bound vectors in MPS input.

\subsection{Validation campaigns}
\label{appE:validation}
\label{app:validation-campaigns}
\paragraph{Targeted invalid certificates and resource limits.}
The targeted suite modifies valid small certificates, with 256 fixtures per
category. All 12,544 fixtures in the 49 invalid categories return
\texttt{Invalid}. Two additional categories impose insufficient work or
stored-statement budgets on otherwise valid certificates; all 512 return
\texttt{ResourceLimit}. Table~\ref{tab:appE:targeted} gives the category groups.
The suite includes out-of-scope statements, omitted destination variables,
invalid arithmetic, false contradiction claims, incompatible reconstruction
premises, and malformed encodings.

\begin{table}[htbp]
\centering
\small
\caption{Targeted validation categories, with 256 fixtures per category.
Resource-limit fixtures contain otherwise valid certificates.}
\label{tab:appE:targeted}
\begin{tabularx}{\linewidth}{@{}>{\raggedright\arraybackslash}Xrrl@{}}
\toprule
Category group & Categories & Fixtures & Observed result \\
\midrule
Core inference, case scope, branching, and complete output arity
 & 18 & 4,608 & All \texttt{Invalid} \\
Local reconstruction rules: premises and dimensions
 & 7 & 1,792 & All \texttt{Invalid} \\
Network reconstruction rules: dimensions, endpoints, costs, and premises
 & 17 & 4,352 & All \texttt{Invalid} \\
Binary encoding: opcodes, rational canonicality, truncation, and trailing bytes
 & 7 & 1,792 & All \texttt{Invalid} \\
Work and stored-statement limits
 & 2 & 512 & All \texttt{ResourceLimit} \\
\bottomrule
\end{tabularx}
\end{table}

\paragraph{Model-pair binding.}
The 240 valid benchmark certificates are cyclically reassigned to incorrect
source--target pairs. All 240 checks return \texttt{Invalid}.

\paragraph{Exact model import.}
The 14 MPS and 10 LP regressions cover objective constants and maximization,
integer markers, default and negative upper bounds, binary-domain intersections,
ranged rows, duplicate coefficients, multiline LP syntax, malformed inputs,
and resource outcomes. All 24 pass. An independent comparison uses 160
generated LP models and canonical models constructed with Python's
\texttt{Fraction} type. All 160 match, with expected rational coefficients
computed independently of the production scalar parser.

\paragraph{Network reconstruction rules.}
The small-network suite contains 128 Flow, 128 Subflow, and 896 UnitFlowRound
instances, with at most five nodes and eight directed edges. All instances pass.
The UnitFlowRound inputs include 257 fractional and 639
integral flows; known binary witnesses are used to check the balance and
cost relations.

\paragraph{Deterministic certificate assembly.}
Assembling each of the 240 RPA certificates twice produces byte-identical RPX
files matching the stored certificates. All 240 pass independent checker
invocations. The default assembler ceiling suffices for 238 inputs;
\texttt{k1mushroom} and \texttt{supportcase12} increase the assembler's syntax
and expression-node budget from 20 million to 64 million. The larger budget
applies only to assembly.

\paragraph{Binary-stream mutation.}
In 50,000 bounded mutation trials starting from three valid certificates, RPX
accepts 429 inputs: 277 unchanged certificates and 152 modified certificates
that remain valid. Independent exact-arithmetic checks validate all 152 modified
certificates by checking the constant inequalities or enumerating both binary
choices for the one-node UnitFlowRound case.

\Needspace{30\baselineskip}
\subsection{Implementation size}
\label{appE:code}
\label{app:trusted-code-audit}
Table~\ref{tab:appE:code} counts nonblank, noncomment physical source lines,
excluding tests. The RPX acceptance implementation comprises 3,893 lines across
the model, certificate, codec, verifier, support, and reconstruction-rule modules.

\begin{table}[H]
\centering
\small
\caption{Implementation size. The reconstruction-rule modules are
included in the acceptance-implementation total.}
\label{tab:appE:code}
\begin{tabularx}{\linewidth}{@{}>{\raggedright\arraybackslash}Xr>{\raggedright\arraybackslash}X@{}}
\toprule
Code region & Source lines & Counting scope \\
\midrule
RPX-owned acceptance implementation & 3,893
 & Complete owned acceptance modules \\
Six reconstruction-rule modules & 1,380
 & Included in the 3,893-line total \\
RPA compiler and assembler & 10,128
 & Authoring software outside the acceptance boundary \\
\bottomrule
\end{tabularx}
\end{table}

The count excludes the CLI and I/O facade, text-import utilities, third-party
arithmetic libraries, and platform software. The artifact documents these dependencies.

\subsection{Reproduction workflow}
\label{appE:reproduction}
\label{app:reproduction-workflow}
The artifact provides the implementation, all 240 model pairs and certificates,
validation suites, raw measurements, and reproduction scripts. A checksum
manifest identifies the distributed files. \texttt{REPRODUCE.txt} lists
reproduction commands; \texttt{docs/FORMAT.md} specifies grammars, file formats,
and resource policies. The numeric header in Listing~\ref{lst:running-rpa}
identifies the serialization format. A minimal workflow from the artifact root is:
\begin{verbatim}
cargo build -p kernel-core --release --locked
target/release/rpx canonical source.mps source.mcm
target/release/rpx canonical target.lp target.mcm
target/release/rpx assemble proof.rpa proof.rpx
target/release/rpx verify source.mcm target.mcm proof.rpx
\end{verbatim}

\end{document}